\documentclass[12pt]{article}
\usepackage{amsmath,amssymb,amsthm}
\usepackage[all]{xy}

\advance\textheight by20mm
\advance\textwidth by20mm
\advance\hoffset by-10mm
\advance\voffset by-10mm

\begin{document}

\newcommand{\GL}{{\mathrm{GL}}}
\newcommand{\SL}{{\mathrm{SL}}}
\newcommand{\Gr}{{\mathrm{Gr}}}
\newcommand{\Fl}{{\mathrm{Fl}}}
\newcommand{\Hom}{\mathop{\mathrm{Hom}}\nolimits}
\newcommand{\Ind}{\mathop{\mathrm{Ind}}\nolimits}
\newcommand{\Int}{\mathop{\mathrm{Int}}\nolimits}
\newcommand{\Init}{\mathop{\mathrm{Init}}\nolimits}
\newcommand{\Ker}{\mathop{\mathrm{Ker}}\nolimits}
\newcommand{\Nuc}{\mathop{\mathrm{Nuc}}\nolimits}
\newcommand{\Stab}{\mathop{\mathrm{Stab}}\nolimits}
\newcommand{\Term}{\mathop{\mathrm{Term}}\nolimits}
\newcommand{\rk}{\mathop{\mathrm{rk}}\nolimits}
\newcommand{\Rep}{{\mathcal{R}\mathrm{ep}}}
\newcommand{\InjRep}{\mathcal{I}{\mathrm{nj}\mathcal{R}\mathrm{ep}}}
\newcommand{\codim}{{\mathrm{codim}}}
\newcommand{\mult}{{\mathrm{mult}}}

\newcommand{\cF}{{\mathcal{F}}}
\newcommand{\cO}{{\mathcal{O}}}
\newcommand{\fA}{{\mathfrak{A}}}
\newcommand{\fB}{{\mathfrak{B}}}
\newcommand{\fC}{{\mathfrak{C}}}
\newcommand{\fF}{{\mathfrak{F}}}

\newcommand{\ba}{{\mathbf{a}}}
\newcommand{\bb}{{\mathbf{b}}}
\newcommand{\bc}{{\mathbf{c}}}
\newcommand{\bd}{{\mathbf{d}}}

\newcommand{\NN}{{\mathbb{N}}}
\newcommand{\KK}{{\mathbb{K}}}
\renewcommand{\AA}{{\mathbb{A}}}
\newcommand{\ssn}{\subsetneq}

\newcommand{\degless}{\stackrel{{\mathrm{deg}}}{\leq}}
\newcommand{\rkless}{\stackrel{{\mathrm{rk}}}{\leq}}
\newcommand{\mvless}{\stackrel{{\mathrm{mv}}}{\leq}}

\newtheorem{theorem}{Theorem}
\newtheorem{prop}[theorem]{Proposition}
\newtheorem{corollary}[theorem]{Corollary}
\newtheorem{lemma}[theorem]{Lemma}

\newcommand{\defin}{{\noindent\bf Definition. }}
\newcommand{\remark}{{\noindent\bf Remark. }}
\newcommand{\example}{{\noindent\bf Example. }}

\title{Bruhat order for two subspaces and a flag}
\author{Evgeny Smirnov}

\maketitle

\abstract{The classical Ehresmann--Bruhat order describes the
possible degenerations of a pair of flags in a finite-dimensional
vector space $V$; or, equivalently, the closure of an orbit of the
group $\GL(V)$ acting on the direct product of two full flag
varieties.

We obtain a similar result for triples consisting of two subspaces and
a partial flag in $V$; this is equivalent to describing the closure of a
$\GL(V)$-orbit in the product of two Grassmannians and one flag
variety. We give a rank criterion to check whether such a triple can
be degenerated to another one, and we classify the minimal degenerations. Our
methods involve only elementary linear algebra and
combinatorics of graphs (originating in Auslander--Reiten quivers).}

\section{Introduction}

We will consider certain configurations of subspaces in an
$n$-dimensional vector space $V$ over an algebraically closed field
$\KK$. These configurations $(U,W,V_\bullet)$ consist of two subspaces
$U$ and $W$ of $V$ of fixed dimensions $k$ and $l$, and a partial flag
$V_\bullet=(V_{d_1}\subset V_{d_2}\subset\dots\subset V_{d_m}=V)$, 
where $\dim V_{d_i}=d_i$.

Our goal is to describe such configurations up to a linear change of
coordinates in $V$ and the ways how configurations
degenerate. In other words, we consider the
direct product $X=\Gr(k,V)\times\Gr(l,V)\times\Fl_{\mathbf d}(V)$ of two
Grassmannians and a flag variety of type $\bd=(d_1,\dots,d_m)$ in $V$, the group $\GL(V)$
acting diagonally on this variety, and describe orbits of this
action and the inclusion relations between their closures.

One can easily show that the number of these orbits is finite. Such a product
$X$ of flag varieties is said to be a \emph{multiple flag variety of finite type}. In the
paper \cite{MWZ} the authors list all such varieties and describe a way of indexing the orbits of the general linear group acting on them.

They also obtain a necessary condition for the closure of a $\GL(V)$-orbit on such a variety to contain another $\GL(V)$-orbit. This condition comes from the results by C.~Riedtmann \cite{R} on degenerations of representations of quivers.

It is not always clear whether this condition provides a criterion. As is mentioned in \cite{MWZ}, this is so in several cases, as follows from some general results on quivers due to K.~Bongartz (\cite[\S 4]{B1}, \cite[\S 5.2]{B2}). One more case is treated in the paper \cite{M} by P.~Magyar, where a similar criterion is obtained for configurations of two flags and a line. Magyar's approach is elementary; it uses only combinatorics and linear algebra.

The case $X=\Gr(k,V)\times\Gr(l,V)\times\Fl_{\bd}(V)$ we are interested in is covered by the results of Bongartz. However, in this case we provide a simpler criterion for a configuration to degenerate to another one, in terms of dimensions of certain subspaces obtained from $U$, $W$, and $V_\bullet$ by taking sums and intersections, and we give a completely elementary proof of this result.


For this, we follow in general the approach of \cite{M}. But the combinatorics we use for indexing the orbits in $X$ is quite different.

For a geometric study of orbit closures in $X$ in the particular case $\bd=(1,\dots,n)$ (that is, when $\Fl_{\bd}(V)$ is the full flag variety; this case we call \emph{spherical}), we address the reader to our paper \cite{Sm}.

\paragraph*{Structure of the paper.} This paper is organised as follows. In Section~\ref{sec2}, we recall some
results from \cite{MWZ} concerning classification of orbits in an
arbitrary multiple flag variety of finite type. In Section~\ref{sec3}, we
introduce an indexing of orbits of $\GL(V)$ in $X$ by subsets of vertices of a certain quiver. Section~\ref{sec4} is devoted to defining three
partial orders on this set of orbits: the first order is
given by degenerations of orbits, the second one is given by
conditions on dimensions of certain subspaces, and the definition of
the third order is purely combinatorial, involving the description of
orbits from Section~\ref{sec3}. In Section~\ref{spherical}, we discuss the relation of the third order with the ``weak order" on spherical varieties in the spherical case. The principal result of this paper states that
the first three orders are the same; this is proved in Section~\ref{sec5}.

\paragraph*{Acknowledgements.} I am grateful to Grzegorz Zwara for extremely useful discussions on Auslander--Reiten quivers, and to Andrei Zelevinsky for drawing my attention to the paper \cite{M}. I also would like to thank Michel Brion for constant attention to this work.

\section{Orbits and representations: a general approach}\label{sec2}

In this section, we consider the problem of classifying orbits of the
general linear group in a multiple flag variety in its general
setting, after \cite{MWZ}.

Let $V$ be an $n$-dimensional vector space over a field $\KK$, which
we suppose to be arbitrary throughout this and the next Section. Let $Q_{p,q,r}$ be the three-arm
star-like quiver of the following form:
$$
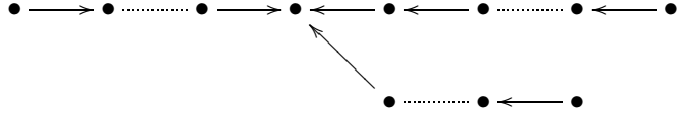
\xymatrix{
{\bullet}\ar[r]&{\bullet}\ar@{.}[r]&{\bullet}\ar[r]&{\bullet}&{\bullet}\ar[l]&{\bullet}\ar[l]\ar@{.}[r]&{\bullet}&{\bullet}\ar[l]\\
&&&&{\bullet}\ar[ul]\ar@{.}[r]&{\bullet}\ar@{.}[r]&{\bullet}\ar[l] }
$$
with $p+q+r-2$ vertices forming three arms of lengths $p$, $q$, and
$r$, and with all arrows leading to the center.

Let $\Rep(Q_{p,q,r})$ denote the category of representations of this
quiver. Magyar, Weyman, and Zelevinsky \cite{MWZ} consider the full
subcategory $\InjRep(Q_{p,q,r})$ in $\Rep(Q_{p,q,r})$ whose objects are those
representations such that all the linear maps corresponding to the
arrows are injections. The subcategory $\InjRep(Q_{p,q,r})$ is closed
under taking direct sums and subobjects (but not quotients!), so one
can introduce the notions of decomposition into direct sums and
indecomposable objects. The uniqueness of a decomposition into a
sum of indecomposables is guaranteed by general results due to Kac
\cite{K}.

In particular, the set of indecomposables $\Ind(\InjRep(Q_{p,q,r}))$
forms a subset of $\Ind(\Rep(Q_{p,q,r}))$, since it is closed under
taking subobjects.

Fix a dimension vector $(\ba,\bb,\bc)=(a_1,\dots,a_p;b_1,\dots,b_q;c_1,\dots,c_r)$, where $a_p=b_q=c_r$, and take a representation 
$$
\underline
V=(V_1,\dots,V_p;V_1',\dots,V_q';V_1'',\dots,
V''_r)\in \InjRep(Q_{p,q,r})
$$
with dimension vector $(\ba,\bb,\bc)$. This representation can be considered as a triple
of partial flags in $V=V_p=V_q'=V_r''$ with the given depths and
dimension vectors, defined up to $\GL(V)$-action. And, vice versa,
any such triple of flags provides a representation from
$\InjRep(Q_{p,q,r})$. So, the orbits of the diagonal action of
$\GL(V)$ on the direct product of three partial flag varieties
$$
\Fl_{(\ba,\bb,\bc)}(V)=\Fl_\ba(V)\times\Fl_\bb(V)\times \Fl_\bc(V)
$$
are in 
one-to-one correspondence with the elements of $\InjRep(Q_{p,q,r})$
with dimension vector $(\ba,\bb,\bc)$.

In this category we  have the uniqueness of a
decomposition into a sum of indecomposables. We also have the following property: there exists at most one indecomposable object  with a given dimension vector. This
means that the $\GL(V)$-orbits in $\Fl_{(\ba,\bb,\bc)}(V)$ correspond to
the possible decompositions of the dimension vector $(\ba,\bb,\bc)$:
$$
(\ba,\bb,\bc)=\sum\underline\dim I_\alpha,
$$
where $I_\alpha$ are indecomposable objects. So, if the number of
$\GL(V)$-orbits in $\Fl_{(\ba,\bb,\bc)}(V)$ is finite (in this case this
multiple flag variety is said to be \emph{of finite type}), the
classification of orbits is thus reduced to a purely combinatorial
problem.

So, knowing all the indecomposable objects in the category
$\InjRep(Q_{p,q,r})$ for a given quiver $Q_{p,q,r}$ allows us to describe the
$\GL(V)$-orbits in the multiple flag variety $\Fl_{(\ba,\bb,\bc)}(V)$ for
an arbitrary dimension vector $(\ba,\bb,\bc)$. The complete list of
all multiple flag varieties of finite type and indecomposable
objects in the corresponding categories is given in
\cite[Theorem~2.3]{MWZ}.

In particular, this list includes quivers $Q_{p,q,1}$ (type $A$) and $Q_{p,2,2}$ (type $D$). The multiple flag varieties corresponding to these two series of quivers will be the main objects of our interest throughout this paper.

\section{Combinatorial enumeration of objects with a specific dimension vector}\label{sec3}

Consider the Auslander--Reiten quiver (AR-quiver) for the category
$\Rep(Q)$. Its vertices correspond to indecomposable objects, and
arrows represent ``minimal'' morphisms between indecomposables
--- i.e., morphisms
$$
f\colon I\to I'
$$
that cannot be presented as a composition of two morphisms
$$
f=g\circ h\colon I\stackrel{h}{\to} I''\stackrel{g}{\to} I',
$$
where $I$, $I'$ and $I''$ are pairwise non-isomorphic
indecomposables.

Having the AR-quiver for $\Rep(Q)$, consider its subquiver
defined as follows: we take all vertices that correspond
to indecomposable objects from $\InjRep(Q)$ and all arrows between these
vertices. This is the Auslander--Reiten quiver for the category
$\InjRep(Q)$. We will refer to \emph{the latter} quiver (not to the
former) as to the AR-quiver for the quiver $Q$; it will be denoted
by $AR(Q)$.

For background on Auslander--Reiten quivers, see the book \cite{ARS}.

Now let us pass to the explicit study of cases $A$ and $D$.

\subsection{Case A: two flags}

Let $Q$ equal $Q_{p,q,1}$. That is, $Q$ is a linear quiver with
$p+q-1$ vertices and arrows oriented as follows:
$\xymatrix@=8pt{{\bullet}\ar[r]&{\bullet}\ar@{.}[r]&{\bullet}\ar[r]&{\stackrel{p}{\bullet}}&
{\bullet}\ar@{.}[r]\ar[l]&{\bullet}&{\bullet}\ar[l]}$

All the indecomposable injective representations of this quiver are
one-di\-men\-sional. They are as follows:
$$\label{indec_A}
I_{ij}=\left(\xymatrix@=8pt{{0}\ar@{.}[r]&{0}\ar[r]&{\KK}\ar@{.}[r]&{\KK}\ar[r]&{\KK}&
{\KK}\ar@{.}[r]\ar[l]&{\KK}&{0}\ar[l]&0\ar@{.}[l]}\right),
$$
where the first nonzero space has number $i$, the last --- the
number $p+q-j$, and $i\in [1,p]$, $j\in [1,q]$. So,
there are $pq$ non-isomorphic indecomposable objects.

The AR-quiver for such a quiver is a rectangle of size $(p\times q)$. Let us
draw the example where $p=4$, $q=3$:

$$
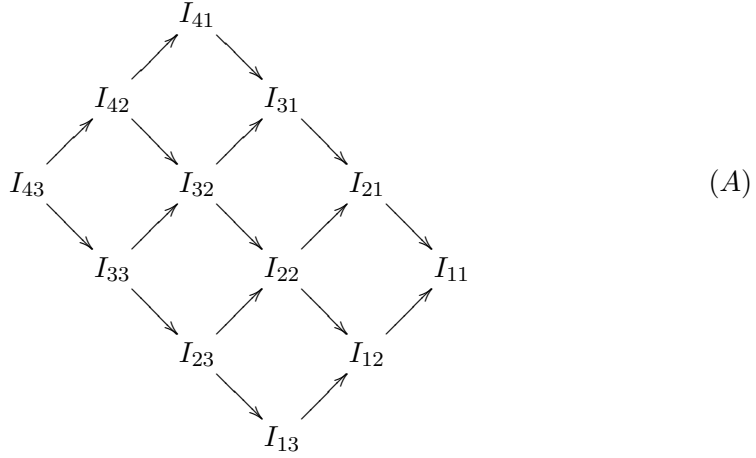
\xymatrix@ur{
{I_{43}}\ar[r]\ar[d]& {I_{42}}\ar[r]\ar[d]& {I_{41}}\ar[d]\\
{I_{33}}\ar[r]\ar[d]& {I_{32}}\ar[r]\ar[d]& {I_{31}}\ar[d]\\
{I_{23}}\ar[r]\ar[d]& {I_{22}}\ar[r]\ar[d]& {I_{21}}\ar[d]\\
{I_{13}}\ar[r]& {I_{12}}\ar[r]& {I_{11}}\\
}\eqno{(A)}
$$

Given an object $F\in\InjRep(Q)$, we will say that an indecomposable object $I$ \emph{occurs} in $F$, if it occurs with nonzero multiplicity in the decomposition of $F$ into indecomposables.

\begin{prop}
Let $F$ be an object in $\InjRep(Q_{p,q,1})$ corresponding to a
configuration of two flags,
such that $\underline\dim F=(a_1,\dots,a_p;b_1,\dots,b_q)$, $a_p=b_q=n$, and let
$F=\bigoplus I_{ij}$ be its decomposition into a sum of
indecomposable objects. Then there are $n$ summands. On each path formed by the elements $I_{i\alpha}$ with $i$ fixed, there are exactly $a_i-a_{i-1}$ indecomposable objects, counted with multiplicities, occuring in $F$. On each path formed by the elements $I_{\alpha j}$ with $j$ fixed, there are exactly $b_j-b_{j-1}$ indecomposable objects occuring in $F$. (We set formally $a_0=b_0=0$).
\end{prop}

\proof Since all the indecomposable summands are one-dimensional,
there are exactly $n$ of them. As we have seen before,
$$
\underline{\dim} I_{ij}=(\underbrace{0,\dots,0,}_{i-1\text{
entry}}1,\dots,1,\dots,1, \underbrace{0,\dots,0,}_{j-1\text{
entry}}).
$$
The resulting dimension is the sum of dimensions of the
indecomposable objects occuring in $F$:
$$
\underline{\dim} F=\sum\underline{\dim} I_{ij}.
$$

Denote the dimension vector of a representation by $(\ba',\bb')=(a_1',\dots,a_p';b_1',\dots,b_q')$. For a given $i$, the objects $I_{ij}$ are characterized by the
equality $a'_{i}=a'_{i-1}+1$. For all other indecomposable objects,
$a'_{i}=a'_{i-1}$. This means that there are
exactly $a_i-a_{i-1}$ objects of the form $I_{ij}$ occuring in $F$.

The fact that $F$ contains exactly $b_j-b_{j-1}$ summands of the form $I_{ij}$ for
a given $j$ is proved similarly.
\endproof

\begin{corollary} Consider the particular case $p=q=n$, $(\ba,\bb)=(1,2\dots,n;1,2\dots,n)$.
Then for any two summands $I_{ij}$ and $I_{i'j'}$ occuring in $F$, we have $i\neq i'$ and
$j\neq j'$. So, objects with such dimension vector 
are in one-to-one correspondence with the configurations of $n$
rooks not attacking each other on the chessboard of size $n\times
n$, i.e., with the permutations of the set of $n$ elements. In
particular, there are $n!$ such non-isomorphic objects.
\end{corollary}

We will see in Section~\ref{spherical} that this description coincides with the well-known indexing of $B$-orbits in a full flag variety by permutations.

\subsection{Case D: two subspaces and a flag}

Now let $Q$ be the quiver $D_{p+2}$ with all arrows mapping to the
center.

Having a representation
$$
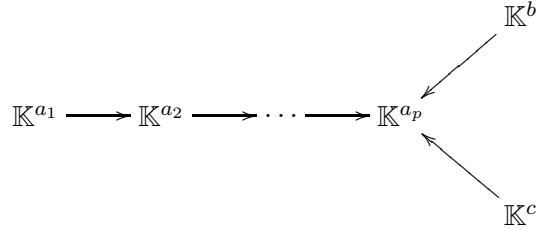
\xymatrix{
&&&&{\KK^{b}}\ar[dl]\\
{\KK^{a_1}}\ar[r]& {\KK^{a_2}}\ar[r]& {\dots}\ar[r]&{\KK^{a_p}}\\
&&&&{\KK^{c}}\ar[ul] }
$$
we denote its dimension vector by $(a_1,\dots,a_p;b;c)$.

Here is the complete list of indecomposable objects in $\InjRep(Q)$,
taken from \cite[Theorem~2.3]{MWZ}.
There are four one-dimensional series, which we present in the table
below together with their dimension vectors:
$$\label{indec_D}
\begin{tabular}{ll}
$I^+_i$ & $(0,\dots,0,1,\dots,1;1;0)$\\
$I^-_i$ & $(0,\dots,0,1,\dots,1;0;1)$\\
$I_{i\infty}$\qquad & $(0,\dots,0,1,\dots,1;0;0)$\\
$I_{0i}$ & $(0,\dots,0,1,\dots,1;1;1)$\\
\end{tabular}
$$
(all the maps between one-dimensional spaces are nonzero, the
dimension jumps at the $i$-th step, $i\in[1,p]$), and one family of
the following form:
$$
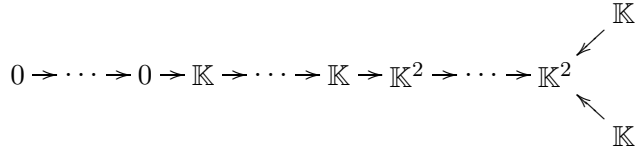
\xymatrix@=9pt{
&&&&&&&&&{\KK}\ar[dl]\\
0\ar[r]&{\dots}\ar[r]&0\ar[r]&{\KK}\ar[r]^{~}&{\dots}\ar[r]&{\KK}\ar[r]&{\KK^2}\ar[r]
&{\dots}\ar[r]&{\KK^2}\\
&&&&&&&&&{\KK}\ar[ul]\\
}
$$
where all the images of the three maps $\KK\to\KK^2$ are distinct
(this guarantees indecomposability), and the dimension within the
longest arm jumps at the $i$-th and the $j$-th steps, $i<j$. Denote
these objects by $I_{ij}$.

From the definition of $AR(Q)$ we obtain the following example, where
$p=5$: 
$$
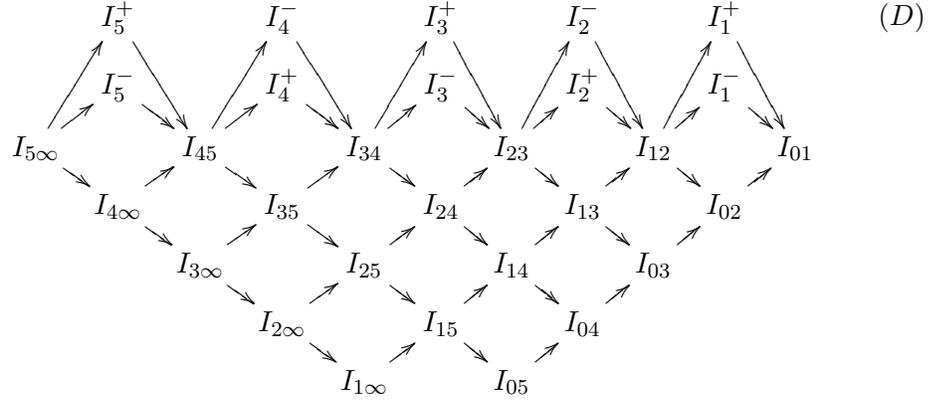
\xymatrix@=7pt{
& {I_5^+}\ar[ddr] &&{I_4^-}\ar[ddr] &&{I_3^+}\ar[ddr] &&{I_2^-}\ar[ddr] &&{I_1^+}\ar[ddr]\\
& {I_5^-}\ar[dr] &&{I_4^+}\ar[dr] &&{I_3^-}\ar[dr] &&{I_2^+}\ar[dr] &&{I_1^-}\ar[dr]\\
{I_{5\infty}}\ar[uur] \ar[ur] \ar[dr] && {I_{45}} \ar[uur] \ar[ur]
\ar[dr] && {I_{34}} \ar[uur] \ar[ur] \ar[dr]
&& {I_{23}}\ar[uur] \ar[ur] \ar[dr] && {I_{12}}\ar[uur] \ar[ur] \ar[dr] && {I_{01}}\\
&{I_{4\infty}}\ar[ur] \ar[dr] && {I_{35}}\ar[ur] \ar[dr] &&
{I_{24}}\ar[ur] \ar[dr] && {I_{13}}\ar[ur] \ar[dr] &&
{I_{02}}\ar[ur]\\
&&{I_{3\infty}}\ar[ur] \ar[dr] && {I_{25}}\ar[ur] \ar[dr] && {I_{14}} \ar[ur] \ar[dr] &&  {I_{03}}\ar[ur]\\
&&&{I_{2\infty}}\ar[ur]\ar[dr] && {I_{15}}\ar[ur]\ar[dr]&&  {I_{04}}\ar[ur]\\
&&&& {I_{1\infty}}\ar[ur]&& I_{05}\ar[ur] }\eqno{(D)}
$$

Indeed, knowing the AR-quiver for $\Rep(D_{p+2})$ with arrows
oriented to the center, we restrict ourselves to its vertices
corresponding to indecomposable objects from
$\InjRep(D_{p+2})$. Construction of the AR-quiver for $\Rep(Q)$ with $Q$ arbitrary is discussed, for instance, in \cite[Chap.~VII]{ARS}

\noindent{\bf Notation. } The two subsets of vertices of the two
top rows connected by the dashed and the dotted line, formed by the
objects of the form $I^+_i$ and $I^-_i$, are called
\emph{zigzags}. Subsets of vertices of the following form, represented
by white circles on the figure below, are said to be
\emph{roads}\label{defroad}:
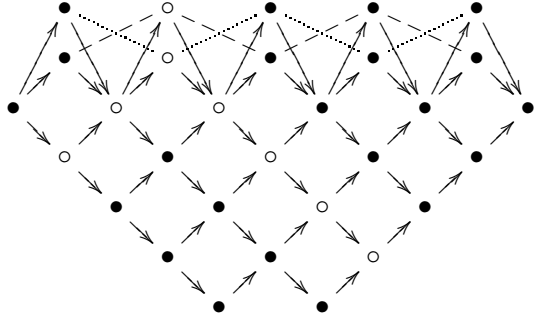
$$
\xymatrix@=8pt{
& {\bullet}\ar[ddr]\ar@{.}[drr] &&{\circ}\ar[ddr]\ar@{--}[drr] &&{\bullet}\ar[ddr]\ar@{.}[drr] &&{\bullet}\ar[ddr]\ar@{--}[drr] &&{\bullet}\ar[ddr]\\
& {\bullet}\ar[dr]\ar@{--}[urr] &&{\circ}\ar[dr]\ar@{.}[urr] &&{\bullet}\ar[dr]\ar@{--}[urr] &&{\bullet}\ar[dr]\ar@{.}[urr] &&{\bullet}\ar[dr]\\
{\bullet}\ar[uur] \ar[ur] \ar[dr] && {\circ} \ar[uur] \ar[ur]
\ar[dr] && {\circ} \ar[uur] \ar[ur] \ar[dr]
&& {\bullet}\ar[uur] \ar[ur] \ar[dr] && {\bullet}\ar[uur] \ar[ur] \ar[dr] && {\bullet}\\
&{\circ}\ar[ur] \ar[dr] && {\bullet}\ar[ur] \ar[dr] &&
{\circ}\ar[ur] \ar[dr]
&& {\bullet}\ar[ur] \ar[dr] &&  {\bullet}\ar[ur]\\
&&{\bullet}\ar[ur]\ar[dr] && {\bullet}\ar[ur]\ar[dr] && {\circ} \ar[ur]\ar[dr] && {\bullet}\ar[ur]\\
&&&{\bullet}\ar[ur]\ar[dr] && {\bullet}\ar[ur]\ar[dr]&& {\circ} \ar[ur]\\
&&&&{\bullet}\ar[ur]&&{\bullet}\ar[ur] }
$$ 
They are formed by the objects $I_{i\infty},\dots,I_{i,i+1}$,
$I^+_i$, $I^-_i$, $I_{i-1,i},\dots,I_{0i}$ for a given $i$.
Each road starts on the left edge of the AR-quiver, at an object $I_{i\infty}$, goes up, then
passes through the ``mountain range" formed by two upper rows,
bifurcates there and then goes down to the right edge, ending at the object $I_{0i}$. This road is said to be the $i$-th one. So, there are exactly 2 different zigzags and $p$ different roads.

\begin{prop} Let $F$ be an object in $\InjRep(Q_{p,2,2})$, such that
$$
\underline{\dim} F=(a_1,a_2,\dots,a_p;k;l),
$$
and let $F=\bigoplus
I_\alpha$ be its decomposition into a sum of indecomposables. Then:
\begin{enumerate}

\item[(i)] For the $i$-th road in $AR(Q_{p,2,2})$ there are exactly $a_i-a_{i-1}$ objects occuring in $F$ situated on this road (as before, $a_0$ is set to be equal to $0$);

\item[(ii)] The total number of $I_\alpha$ of the form $I_{ij}$, $1\le i<j\le
  n$, and $I^+_i$, equals $k$;

\item[(iii)] The total number of $I_\alpha$ of the form $I_{ij}$, $1\le i<j\le
  n$, and $I^-_i$, equals $l$.

\end{enumerate}
\end{prop}

\proof Fix a road; let $I_{i\infty}$ be its first element. From the
description of indecomposable objects given on 
Page~\pageref{indec_D}, it follows that the dimension vectors $(\ba';b';c')$ of the
indecomposable objects situated on this road are characterized by the
equality $a'_i=a'_{i-1}+1$. For all other elements, $a'_i=a'_{i-1}$. So, $F$ contains exactly $a_i-a_{i-1}$ indecomposable objects with dimension jump on the $i$-th step. This
proves the first part of the proposition.

(ii) and (iii) are proved similarly.\endproof

So, an object with  dimension vector $(a_1,\dots,a_p;k;l)$ gives us a
set of vertices in $AR(D_{p+2})$, satisfying the properties
(i)--(iii). Obviously, the converse is also true: each set of
vertices determines an object, namely, the direct sum of the
corresponding indecomposables, and the properties (i)--(iii)
guarantee that the dimension vector of this object equals
$(a_1,\dots,a_p;k;l)$.\endproof

\section{Three orders}\label{sec4}

Throughout this section, $Q$ is either the quiver
$A_{p+q-1}=Q_{p,q,1}$ or the quiver $D_{p+2}=Q_{p,2,2}$. Recall that throughout the rest of this paper, the ground field $\KK$ is supposed to be algebraically closed.

In this section we present three different ways to turn the set of
objects $F\in\InjRep(Q)$ with a given dimension vector into
a partially ordered set (or shortly \emph{poset}). We will show that these three orders are the same in the next section. 


\subsection{Degeneration order}

The first definition uses the bijection between objects with dimension
vector $(\ba,\bb,\bc)$ and orbits in the corresponding multiple flag
variety $\Fl_{(\ba,\bb,\bc)}(V)$. Given an object $F$, we denote the corresponding orbit
by $\cO_F$.

\defin We say that $F$ is less or equal than $F'$ w.r.t. the \emph{degeneration
order}, if there is an inclusion of the corresponding orbit closures
(in the Zariski topology):
$$
F\degless F'\quad \Leftrightarrow \quad\cO_F\subseteq\bar\cO_{F'}.
$$

\subsection{Rank order}

Another partial order is defined by means of dimensions of the
homomorphism spaces between objects in the category $\InjRep(Q)$.
For short, for two elements $F,G\in\InjRep(Q)$ we denote the
dimension $\dim\Hom(F,G)$ by $\langle F,G\rangle$.

\defin $F$ is less or equal than $F'$ w.r.t. the \emph{rank
order} (notation: $F\rkless F'$), if for each indecomposable object
$I\in\InjRep(Q)$
$$
\langle I,F\rangle\geq\langle I,F'\rangle.
$$
(NB: the inequality is reversed!)

In our cases ($A_{p+q-1}$ and $D_{p+2}$) we shall give a simple geometric interpretation of the numbers $\langle I,F\rangle$. In general, this interpretation
also exists (see \cite[Prop.~4.1]{MWZ}), but it is not evident at
all.

\begin{prop}\label{ranks}
\begin{enumerate}
\item Let $Q$ equal $Q_{p,q,1}$, and let $V_\bullet=(V_{a_1}\subseteq\dots\subseteq V_{a_p}=V)$
and $V'_\bullet=(V'_{b_1}\subseteq\dots\subseteq V'_{b_q}=V)$ be two flags
of the same depth in a vector space $V$. Then for the object $F$
corresponding to the configuration $(V_\bullet,V'_\bullet)$ the
following equalities hold:
$$
\langle I_{ij},F\rangle=\dim V_{a_i}\cap V_{b_j}'
$$
for each $i\in[1,p]$, $j\in[1,q]$. (A description of the $I_{ij}$ is given on Page~\pageref{indec_A}.)

\item Let $Q$ equal $Q_{p,2,2}$, and let $V_\bullet=(V_{a_1}\subseteq\dots\subseteq
V_{a_p}=V)$, $U$ and $W$ be a flag and two subspaces in $V$. Then for the
object $F$ corresponding to the configuration $(U,W,V_\bullet)$ the
following equalities hold:

\begin{eqnarray}
\langle I_{i\infty},F\rangle&=&\dim V_{a_i}=a_i;\nonumber\\
\langle I_i^+,F\rangle&=&\dim V_{a_i}\cap U;\nonumber\\
\langle I_i^-,F\rangle&=&\dim V_{a_i}\cap W;\label{expranks}\\
\langle I_{0i},F\rangle&=&\dim V_{a_i}\cap U\cap W;\nonumber\\
\langle I_{ij},F\rangle&=&\dim V_{a_j}\cap U\cap W +\dim V_{a_i}\cap((V_{a_j}\cap U)+(V_{a_j}\cap W)).\nonumber
\end{eqnarray}
\end{enumerate}
\end{prop}

\proof A first observation: these formulas are additive under taking
direct sums of objects and componentwise direct sums of
corresponding configurations of subspaces. 

Next, the bracket $\langle \cdot,\cdot\rangle$ is bilinear, so
$$
\langle I,F\oplus F'\rangle=\langle I,F\rangle+\langle I,F'\rangle.
$$
 Thus, it only suffices to prove these formulas for an indecomposable $F$. And
this is done by a direct verification.\endproof

\defin The numbers $\langle I,F\rangle$ are called \emph{rank
numbers}.

\subsection{Move order}\label{moves}

In the previous section we have obtained a combinatorial description
of objects in $\InjRep(Q)$ with a given dimension vector. Objects
are encoded by subsets of vertices of a certain quiver, satisfying a
number of properties.

To introduce the third partial order, we define some operations,
called \emph{elementary moves}, that bring these subsets of vertices
into other ones.

As usual, we begin with type $A$. In this case the definition of elementary move is quite simple.

Take the decomposition of $F$ into indecomposables: $F=\bigoplus
I_\alpha$. 
Suppose that among these $I_\alpha$'s there are two objects $I_{ij}$
and $I_{i'j'}$ occuring in $F$ (probably with multiplicities), such that $i>i'$ and $j>j'$. Let us also suppose
that there is no other $I_{i''j''}$, such that $i>i''>i'$ and
$j>j''>j'$. Graphically, this can be reformulated as follows: there
is no other vertex occuring in $F$ and situated in the following rectangle:
$$
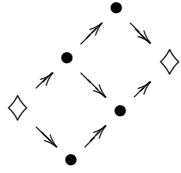
\xymatrix@ur @=8pt{
{\diamondsuit}\ar[r]\ar[d]&{\bullet}\ar[r]\ar[d]&{\bullet}\ar[d]\\
{\bullet}\ar[r]&{\bullet}\ar[r]&{\diamondsuit} }
$$
If this is the case, this rectangle is called \emph{admissible}.

Having this, we construct an object $F'$ by replacing this pair of
indecomposables $I_{ij}\oplus I_{i'j'}$ with the pair $I_{ij'}\oplus
I_{i'j}$. This means that the multiplicities $\mult_{F'}I$ of occurences of  indecomposable objects $I$ in $F'$ are obtained from  $\mult_F I$ according to the following rule:
\begin{eqnarray*}
\mult_{F'} I_{ij}&=&\mult_{F} I_{ij}-1;\\
\mult_{F'} I_{i'j'}&=&\mult_{F} I_{i'j'}-1;\\
\mult_{F'} I_{i'j}&=&\mult_{F} I_{i'j}+1;\\
\mult_{F'} I_{ij'}&=&\mult_{F} I_{ij'}-1;\\
\mult_{F'} I &= &\mult_{F} I\qquad \text{otherwise.}
\end{eqnarray*}

Informally,  can be described as flipping the rectangle, whose ``corners" $I_{ij}$ and $I_{i'j'}$ occuring in $F$ are replaced by $I_{i'j}$ and $I_{ij'}$:
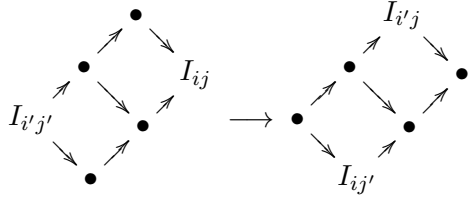
$$
\xymatrix@ur @=8pt{
{I_{i'j'}}\ar[r]\ar[d]&{\bullet}\ar[r]\ar[d]&{\bullet}\ar[d]\\
{\bullet}\ar[r]&{\bullet}\ar[r]&{I_{ij}} }\longrightarrow {
\xymatrix@ur @=8pt{
{\bullet}\ar[r]\ar[d]&{\bullet}\ar[r]\ar[d]&{I_{i'j}}\ar[d]\\
{I_{ij'}}\ar[r]&{\bullet}\ar[r]&{\bullet} }}
$$

Let $F'$ be obtained from $F$ by an elementary move. We denote this
as follows: $F\lessdot F'$.

Now we are ready to give the definition of the third order.

\defin\label{defmoveord} An object $F$ is said to be less or equal than an object
$F'$ w.r.t. the \emph{move order}, if there exists a sequence of
objects $F_0, F_1,\dots, F_s$, such that
$$
F=F_0\lessdot F_1\lessdot\dots\lessdot F_s=F'.
$$
This is denoted as follows: $F\mvless F'$.

\remark Of course, each element is less or equal than itself. This
corresponds to the empty sequence.

So, given two vertices of the AR-quiver, we have at most one
possibility to perform an elementary move affecting them. As a
result of this move, this pair of vertices is replaced with another
pair.

In type $D$ everything is more complicated. As above,
elementary moves consist in replacing a pair of marked vertices, but
now they can be replaced by one, two or three other vertices.
Moreover, the choice of an initial pair does not uniquely define the
move any more; there may be up to three different possibilities.

To begin with, we introduce some convention that allows us to make the
description of elementary moves less bulky. Let us
add a ``fake vertex" in the missing lowest corner, and the
corresponding fake indecomposable object $I_{0\infty}$, equal to
zero. So, the resulting quiver will be as follows:
$$
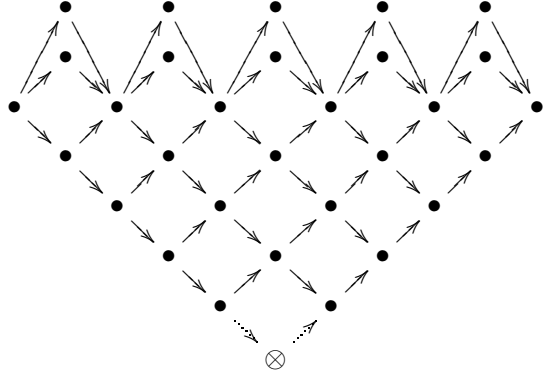
\xymatrix@=8pt{
& {\bullet}\ar[ddr]&&{\bullet}\ar[ddr] &&{\bullet}\ar[ddr]&&{\bullet}\ar[ddr] &&{\bullet}\ar[ddr]\\
& {\bullet}\ar[dr] &&{\bullet}\ar[dr] &&{\bullet}\ar[dr] &&{\bullet}\ar[dr] &&{\bullet}\ar[dr]\\
{\bullet}\ar[uur] \ar[ur] \ar[dr] && {\bullet} \ar[uur] \ar[ur]
\ar[dr] && {\bullet} \ar[uur] \ar[ur] \ar[dr]
&& {\bullet}\ar[uur] \ar[ur] \ar[dr] && {\bullet}\ar[uur] \ar[ur] \ar[dr] && {\bullet}\\
&{\bullet}\ar[ur] \ar[dr] && {\bullet}\ar[ur] \ar[dr] &&
{\bullet}\ar[ur] \ar[dr]
&& {\bullet}\ar[ur] \ar[dr] &&  {\bullet}\ar[ur]\\
&&{\bullet}\ar[ur]\ar[dr] && {\bullet}\ar[ur]\ar[dr] && {\bullet} \ar[ur]\ar[dr] && {\bullet}\ar[ur]\\
&&&{\bullet}\ar[ur]\ar[dr] && {\bullet}\ar[ur]\ar[dr]&& {\bullet} \ar[ur]\\
&&&&{\bullet}\ar[ur]\ar@{.>}[dr]&&{\bullet}\ar[ur] \\
&&&&&{\otimes}\ar@{.>}[ur] }
$$

Now let us describe the moves explicitly. 

Our general strategy will be as follows: first, we define \emph{regions}, which are analogues of rectangles in the case $A_n$. 

A \emph{region} is a triple $(\fA,\Init\fA,\Term\fA)$, where $\fA$ is a subquiver in our AR-quiver of a certain form, described below. Each $\fA$ has exactly one source (vertex of incoming degree 0) and one sink (vertex of outcoming degree 0). These two vertices are called \emph{initial vertices}; we denote this two-elementary set by $\Init\fA$. There are also at least one and at most three vertices marked as \emph{terminal} ones, denoted $\Term\fA$ (they will be defined below in an \emph{ad hoc} way).

\remark The uniqueness of a source and a sink implies, in particular, that $\fA$ is connected and that there exists an (oriented) path joining the initial vertices.


Now let us describe regions explicitly. We distinguish between the following six cases, denoted I.a)-I.e) and II.

The cases I.a)--I.e) are characterized by the following property: $\fA$ consists of those vertices that are situated on the paths joining the source of $\fA$ with its sink.

{\bf I.a)} The initial vertices of a region of type I.a) are of the form $I_1=I_{i'j'}$, $I_2=I_{ij}$, where $i<i'<j<j'$. In this case we define an admissible region $\fA$ of type I.a) as follows:
$$
\fA=\{I_{\alpha\beta}\mid i\leq\alpha\leq i',j\leq\beta\leq j'\}. 
$$
It is a rectangle with corners in $I_1$ and $I_2$. We define the terminal vertices as the two other corners of this rectangle, $I_{ij'}$ and $I_{i'j}$:

A region of this type is shown on the figure. The initial vertices are outlined by  squares, the terminal ones --- by circles.
$$
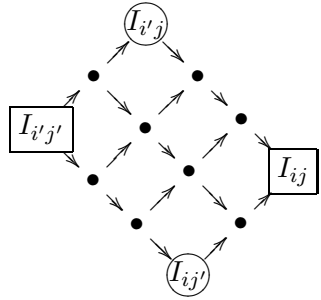
\xymatrix@ur @=8pt{
*+[F]{I_{i'j'}}\ar[r]\ar[d]&{\bullet}\ar[r]\ar[d]&*+[o][F]{I_{i'j}}\ar[d]\\
{\bullet}\ar[r]\ar[d]&{\bullet}\ar[r]\ar[d]&{\bullet}\ar[d]\\
{\bullet}\ar[r]\ar[d]&{\bullet}\ar[r]\ar[d]&{\bullet}\ar[d]\\
*+[o][F]{I_{ij'}}\ar[r]&{\bullet}\ar[r]&*+[F]{I_{ij}} }
$$

{\bf I.b)} The initial vertices of regions of this type are of form $I_1=I_{i'j'}$, $I_2=I_{ij}$, such that $0\leq i<j\leq i'<j'\leq\infty$. For each such pair of vertices, there are two regions of type I.b), defined as follows:
$$
\fA^+=\fA^-=\{I_{\alpha\beta}\mid i\leq \alpha\leq i',j\leq\beta\leq j'\}\cup\{I^+_\gamma,I^-_\gamma \mid j\leq\gamma\leq i'\}
$$

Each such region has three terminal vertices, defined by
\begin{eqnarray*}
\Term\fA^+&=&\{I_{ij'},I_{j}^+,I_{i'}^-\};\\
\Term\fA^-&=&\{I_{ij'},I_{j}^-,I_{i'}^+\}.
\end{eqnarray*}

These two regions are shown on the figures below. 
$$
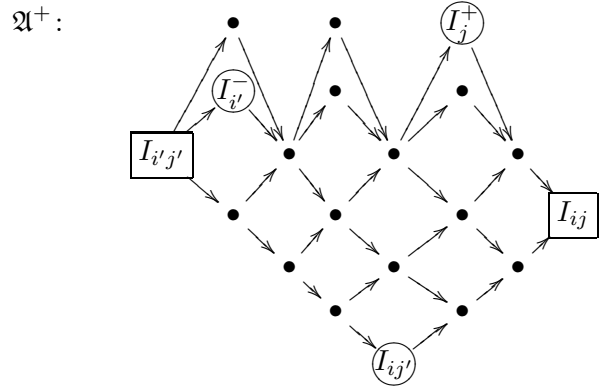
\fA^+\colon\qquad
\xymatrix@=6pt{
& {\bullet}\ar[ddr] &&{\bullet}\ar[ddr] &&*+[o][F]{I^+_{j}}\ar[ddr]\\
& *+[o][F]{I^-_{i'}}\ar[dr] &&{\bullet}\ar[dr] &&{\bullet}\ar[dr]\\
*+[F]{I_{i'j'}}\ar[uur] \ar[ur] \ar[dr] && {\bullet} \ar[uur]
\ar[ur] \ar[dr] && {\bullet} \ar[uur] \ar[ur] \ar[dr]
&& {\bullet}\ar[dr]\\
&{\bullet}\ar[ur] \ar[dr] && {\bullet}\ar[ur] \ar[dr] &&
{\bullet}\ar[ur] \ar[dr]
&& *+[F]{I_{ij}}\\
&&{\bullet}\ar[ur]\ar[dr] && {\bullet}\ar[ur]\ar[dr] && {\bullet} \ar[ur]\\
&&&{\bullet}\ar[ur]\ar[dr] && {\bullet}\ar[ur]\\
&&&&*+[o][F]{I_{ij'}}\ar[ur] }
$$

$$
\fA^-\colon\qquad
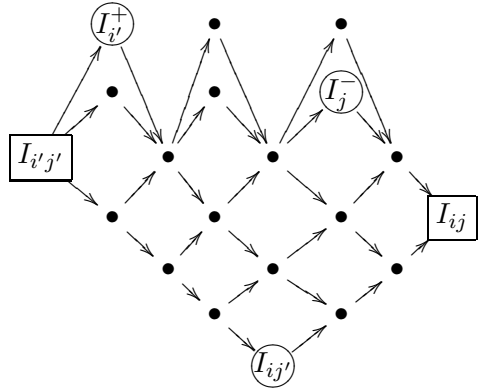
\xymatrix@=6pt{
& *+[o][F]{I^+_{i'}}\ar[ddr] && {\bullet}\ar[ddr] &&{\bullet}\ar[ddr]\\
& {\bullet}\ar[dr] &&{\bullet}\ar[dr] && *+[o][F]{I^-_{j}}\ar[dr]\\
*+[F]{I_{i'j'}}\ar[uur] \ar[ur] \ar[dr] && {\bullet} \ar[uur]
\ar[ur] \ar[dr] && {\bullet} \ar[uur] \ar[ur] \ar[dr]&& {\bullet}\ar[dr]\\
&{\bullet}\ar[ur] \ar[dr] && {\bullet}\ar[ur] \ar[dr] &&
{\bullet}\ar[ur] \ar[dr]
&& *+[F]{I_{ij}}\\
&&{\bullet}\ar[ur]\ar[dr] && {\bullet}\ar[ur]\ar[dr] && {\bullet} \ar[ur]\\
&&&{\bullet}\ar[ur]\ar[dr] && {\bullet}\ar[ur]\\
&&&&*+[o][F]{I_{ij'}}\ar[ur] }
$$

{\bf I.c)} For regions of this type, the initial vertices are of the form $I_1=I_{i'j'}$, $I_2=I^\pm_{i}$, such that $i<i'<j'$. In this case, 
we define $\fA$ to be
$$
\fA=\{I_{\alpha\beta}\mid i\leq \alpha\leq i',\beta\leq j'\}\cup\{I^+_\gamma,I^-_\gamma\mid i\leq\gamma\leq i'\}\cup\{I_{i'}^\pm\},
$$
and $\Term\fA=\{I_{i'}^\pm,I_{ij'}\}$.

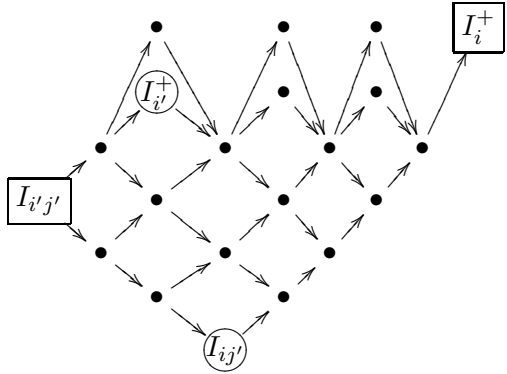
$$
\xymatrix@=6pt{
& &&{\bullet}\ar[ddr] &&{\bullet}\ar[ddr] &&{\bullet}\ar[ddr] &&*+[F]{I^+_{i}}\\
& &&*+[o][F]{I^+_{i'}}\ar[dr] &&{\bullet}\ar[dr] &&{\bullet}\ar[dr]\\
 && {\bullet} \ar[uur] \ar[ur]
\ar[dr] && {\bullet} \ar[uur] \ar[ur] \ar[dr]
&& {\bullet}\ar[uur] \ar[ur] \ar[dr] && {\bullet}\ar[uur]\\
&*+[F]{I_{i'j'}}\ar[ur] \ar[dr] && {\bullet}\ar[ur] \ar[dr] &&
{\bullet}\ar[ur] \ar[dr]
&& {\bullet}\ar[ur]\\
&&{\bullet}\ar[ur]\ar[dr] && {\bullet}\ar[ur]\ar[dr] && {\bullet} \ar[ur]\\
&&&{\bullet}\ar[ur]\ar[dr] && {\bullet}\ar[ur]\\
&&&&*+[o][F]{I_{ij'}}\ar[ur]}
$$

{\bf I.d)} The initial vertices are of the form $I_1=I^\pm_{j'}$, $I_2=I_{ij}$, and $i<j<j'$. Then
$$
\fA=\{I_{\alpha\beta}\mid i\leq \alpha,j<\beta\leq j'\}\cup\{I^+_\gamma,I^-_\gamma \mid j\leq\gamma\leq j'\}\cup\{I_{j'}^\pm\},
$$
and $\Term\fA=\{I_j^\pm,I_{ij'}\}$.
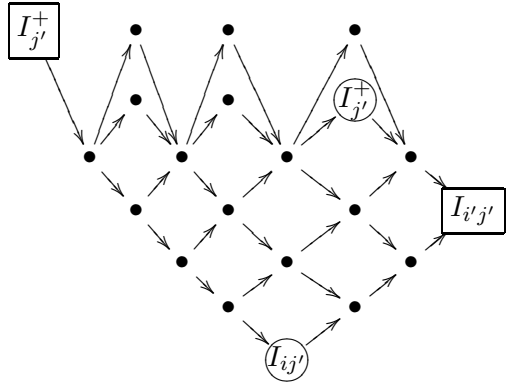
$$
\xymatrix@=6pt{
*+[F]{I^+_{j'}}\ar[ddr] &&{\bullet}\ar[ddr] &&{\bullet}\ar[ddr] &&{\bullet}\ar[ddr]\\
&&{\bullet}\ar[dr] &&{\bullet}\ar[dr] &&*+[o][F]{I^+_{j'}}\ar[dr]\\
& {\bullet} \ar[uur] \ar[ur] \ar[dr] && {\bullet} \ar[uur] \ar[ur] \ar[dr]
&& {\bullet}\ar[uur] \ar[ur] \ar[dr] && {\bullet}\ar[dr]\\
&&{\bullet}\ar[ur] \ar[dr] && {\bullet}\ar[ur] \ar[dr] &&
{\bullet}\ar[ur] \ar[dr] && *+[F]{I_{i'j'}}\\
&&&{\bullet}\ar[ur]\ar[dr] && {\bullet}\ar[ur]\ar[dr] && {\bullet} \ar[ur]\\
&&&&{\bullet}\ar[ur]\ar[dr] && {\bullet}\ar[ur]\\
&&&&&*+[o][F]{I_{ij'}}\ar[ur]}
$$

{\bf I.e)} The initial vertices are of the form $I_{i}^\pm$ and $I_{i'}^\mp$ (signs are different), $i< i'$. Then 
$$
\fA=\{I_{\alpha\beta}\mid i\leq\alpha<\beta\leq i'\}\cup\{I^+_\gamma,I^-_\gamma\mid i<\gamma<i'\}\cup\{I_{i}^\pm,I_{i'}^\mp\}.
$$
Then there is a unique terminal vertex: $\Term\fA=\{I_{ii'}\}$.
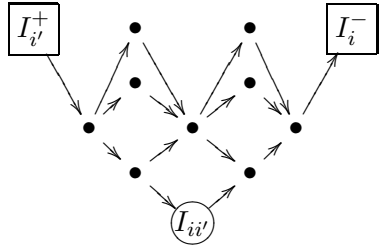
$$
\xymatrix@=6pt{
*+[F]{I^+_{i'}}\ar[ddr]&&{\bullet}\ar[ddr]&&{\bullet}\ar[ddr]&&*+[F]{I^-_i}\\
&&{\bullet}\ar[dr]&&{\bullet}\ar[dr]\\
&{\bullet}\ar[uur]\ar[ur]\ar[dr]&&{\bullet}\ar[uur]\ar[ur]\ar[dr]&&{\bullet}\ar[uur]\\
&&{\bullet}\ar[ur]\ar[dr]&&{\bullet}\ar[ur]\\
&&&*+[o][F]{I_{ii'}}\ar[ur]
}
$$

{\bf II.} In this case, the initial vertices are of the form $I_{ij}$ and $I_{i'j'}$, where $i<j<i'<j'$. The corresponding subquiver $\fA$ is given by
$$
\fA=\{I_{\alpha\beta}\mid i\leq \alpha\leq i',j\leq\beta\leq j'\}\cup\{I^+_\gamma,I^-_\gamma\mid j\leq \gamma\leq i'\},
$$
$I_{i'j}$ and $I_{ij'}$ are its terminal vertices:
$$
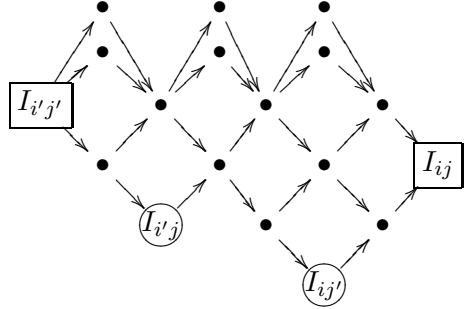
\xymatrix@=6pt{
& {\bullet}\ar[ddr] &&{\bullet}\ar[ddr] &&{\bullet}\ar[ddr]\\
& {\bullet}\ar[dr] &&{\bullet}\ar[dr] &&{\bullet}\ar[dr]\\
*+[F]{I_{i'j'}}\ar[uur] \ar[ur] \ar[dr] && {\bullet} \ar[uur] \ar[ur]
\ar[dr] && {\bullet} \ar[uur] \ar[ur] \ar[dr]
&& {\bullet}\ar[dr]\\
&{\bullet}\ar[ur] \ar[dr] && {\bullet}\ar[ur] \ar[dr] &&
{\bullet}\ar[ur] \ar[dr]
&& *+[F]{I_{ij}}\\
&&*+[o][F]{I_{i'j}}\ar[ur] && {\bullet}\ar[ur]\ar[dr] && {\bullet} \ar[ur]\\
&&&&& *+[o][F]{I_{ij'}}\ar[ur]\\
}
$$
One can think of the obtained set of vertices as a ``folded
rectangle'', with corners in the initial and the terminal vertices.

After having defined regions, we can go further and pass to the definition of the move order. For the following definition, we fix an object $F\in\InjRep(Q_{p,2,2})$.

\defin A region $\fA$ is called \emph{admissible} w.r.t. an object $F$, if for both initial vertices of $\fA$, the corresponding indecomposable objects occur in $F$ with nonzero multiplicities. An admissible region $\fA$ is called \emph{minimal}, if any non-initial vertex from $\fA$ occurs in $F$ with multiplicity 0.

As in the case $A$, elementary moves that can be performed with an object $F$ correspond to the minimal admissible regions:

\defin We say that $F'$ is obtained from $F$ by an \emph{elementary move} (notation: $F\lessdot F'$, if there is a minimal admissible region $\fA$ w.r.t. $F$, such that
\begin{eqnarray*}
\mult_{F'} I =\mult_{F} I-1 & \text{for $I\in\Init\fA$};\\
\mult_{F'} I = \mult_{F} I+1 & \text{for $I\in\Term\fA$};\\
\mult_{F'} I = \mult_{F} I & \text{otherwise.}
\end{eqnarray*}
This means that, as a result of an elementary move, a pair of indecomposable objects is replaced by one, two or three other indecomposable objects.

Now the \emph{move order} is defined as follows: $F$ is said to be less or equal than $F'$ (notation: $F\mvless F'$), if $F'$ is obtained from $F$ by a sequence of elementary moves.

\section{The spherical case, $B$-orbits in $\Gr(k,V)\times\Gr(l,V)$, and weak order}\label{spherical}

Throughout this section, we let the dimension vector $\ba$ be $(1,2,\dots,n)$, so $\Fl_{\ba}(V)$ equals the full flag variety $\Fl(V)$.

Instead of studying orbits of $\GL(V)$ acting on $X=\Gr(k,V)\times \Gr(l,V)\times \Fl(V)$, one can consider the stabilizer $B\subset\GL(V)$ of a complete flag $V_\bullet\in\Fl(V)$ (so that $B$ is a Borel subgroup of $\GL(V)$), and the orbits of $B$ acting diagonally on $Y=\Gr(k,V)\times \Gr(l,V)$. There is an evident bijection between these two sets of orbits, that also respects the degeneration order. So, $Y$ is a $\GL(V)$-variety containing finitely many $B$-orbits. For an orbit $\cO$ in $X$, denote by $\cO_Y$ the corresponding orbit in $Y$.

Consider an arbitrary $\GL(V)$-variety $Z$ with a finite number of $B$-orbits (for an arbitrary connected reductive algebraic group $G$, such varieties are called \emph{spherical}). The set of its orbits admits, along with the usual degeneration order given by
$$
\cO_1\degless\cO_2\Leftrightarrow \cO_1\subset\bar\cO_2,
$$
another partial order structure, called the \emph{weak order}. It was first introduced in \cite{RS} for symmetric spaces, and in \cite{Kn} for spherical varieties.

For its definition, we shall use the \emph{minimal parabolic subgroups} in $\GL(V)$, that is, minimal subgroups containing $B$. There are $n-1$ of them; they are of the form
$$
P_i=\Stab_{\GL(V)}V_\bullet^{(i)},
$$
where $V_\bullet^{(i)}$ is the partial flag $V_1\subset\dots\subset V_{i-1}\subset V_{i+1}\subset V_n=V$, obtained from the standard flag $V_\bullet$ by omitting the $i$-th term.

It is interesting to know when the closure of an orbit in $Y$ is obtained from another orbit closure by the action of a minimal parabolic subgroup:
\begin{equation}\label{onestepparab}
\overline{\cO'_Y}=P_i\cdot\overline{\cO_Y}.
\end{equation}
(we suppose that $\overline{\cO'_Y}\neq \overline{\cO_Y}$; in this case $\dim\cO_Y'=\dim\cO_Y+1$).

The following proposition shows that this relation corresponds to elementary moves with certain properties.

\begin{prop}\label{parab} The equality (\ref{onestepparab}) holds iff  for the objects $F$ and $F'$, corresponding to $\cO_Y$ and $\cO'_Y$, 
$$
F\lessdot F',
$$
and, moreover, the corresponding elementary move is of type I.a), I.c), I.d), I.e), or II, and the source and the sink of the corresponding admissible region belong to neighbor roads.
\end{prop}

The proof of this proposition will be given at the end of Subsection~\ref{ssmvdeg}.

Now let us pass to the definition of the weak order. It is similar to the move order, but its ``elementary moves" are given by the relation (\ref{onestepparab}). Namely, $\cO_Y$ is said to be less or equal than $\cO'_Y$, if there exists a sequence $(P_{i_1},\dots,P_{i_r})$ of minimal parabolic subgroups (possibly with repetitions), such that
$$
\bar\cO'_Y=P_{i_r}\dots P_{i_1}\bar\cO_Y.
$$
We denote this as follows: $\cO_Y\preceq\cO'_Y$.

Obviously, if $\cO_Y\preceq\cO'_Y$, then $\cO_Y\degless\cO'_Y$ (this explains the term ``weak"). However, for an arbitrary spherical variety $Z$, the converse is not true. 
For example, the degeneration order admits a unique maximal element, namely, the open $B$-orbit, and the weak order admits a maximal element for each $G$-orbit on $Z$:  the maximal elements for the weak order are those $B$-orbits that are open in the corresponding $G$-orbit\footnote{In general, this is also false for $G$-homogeneous varieties; an example is provided, for instance, by a full flag variety $\Fl(V)$, where $\dim V\geq 3$.}. 
In particular, $Y=\Gr(k,V)\times\Gr(l,V)$ is not $\GL(V)$-homogeneous, so in this case the weak order is strictly weaker than the degeneration one.

In our paper \cite{Sm}, we describe the weak order on the set of $B$-orbits in $Y$ and then use this description for constructing desingularizations of their closures.

\section{The main result}\label{sec5}

\begin{theorem} Let $Q$ equal $Q_{p,2,2}$. Then for all $F,F'\in\InjRep(Q)$, such that
$\underline\dim F= \underline\dim F'$,
$$
F\degless F'\quad\Leftrightarrow\quad F\rkless
F'\quad\Leftrightarrow\quad F\mvless F'.
$$
So, all the three orders are the same.
\end{theorem}

This is proved in \cite{M} for $Q=Q_{p,q,1}$. We follow the same strategy and split the proof into three lemmas, corresponding to \cite[Lemmas 5,6,7]{M}.

\begin{lemma}\label{mvdeg} $F\mvless F' \Longrightarrow F\degless F'$.
\end{lemma}

This will be proved in \ref{ssmvdeg} by constructing an explicit degeneration of the larger
of the corresponding orbits to the smaller one.

\begin{lemma}\label{degrk} $F\degless F' \Longrightarrow F\rkless F'$.
\end{lemma}

This is a particular case of \cite[Prop.~2.1]{R}. However, in~\ref{ssdegrk} we present an elementary geometric proof of this result.

\begin{lemma}\label{rkmv} $F\rkless F'\Longrightarrow F\mvless F'$.
\end{lemma}

This will be proved in~\ref{ssrkmv} as follows: given $F\rkless F'$, we find an object
$\tilde F$, such that $F\stackrel{\mathrm{mv}}{\lessdot}\tilde
F\rkless F'$.

\subsection{Move order implies degeneration order}\label{ssmvdeg}

First let us recall the description of ``standard'' representatives in
$GL(V)$-orbits, taken from \cite[Def.~2.8, Prop.~2.9]{MWZ}. As usual,
this is described on orbits $\cO_I$ corresponding to indecomposable
objects $I$, and then extended via taking direct sums. 

Let $(U,W,V_\bullet)$ be a triple corresponding to an indecomposable
object. This means that $V=V_{a_p}$ is of dimension 1 or 2. If $\dim
V=1$, each of $U$ and $W$ is either equal to $V$ or to
zero. 

If $I=I_{ij}$, $0<i<j<\infty$, then  $\dim V=2$. Let $(e_i,e_j)$ be
an ordered basis of $V$, such that $V_i=\dots=V_{j-1}=\langle
e_i\rangle$. Then the triple $(U,W,V_\bullet)$ with
$U=\langle e_j\rangle$, $W=\langle e_i+e_j\rangle$ is called the standard representative of the orbit $\cO_{I_{ij}}$.

Later on, we will deal with certain deformations of bases in our
subspaces. For this, the following notational convention will be
useful. Introduce two more ``vectors": $e_0$ and $e_\infty$. Set
formally $e_0=0$ and each linear combination of vectors involving
$e_\infty$ be also equal to 0. Note that with this convention, the
definition of standard representatives for $I_{ij}$, $0<i<j<\infty$,
is extended to the cases of $I_{0i}$ and $I_{i\infty}$, so later
we will consider these three cases simultaneously.

Now we pass to the proof of Lemma~\ref{mvdeg}.

\proof[Proof of Lemma \ref{mvdeg}] The main idea is as
follows: for any two objects $F$ and $F'$, such that $F\lessdot F'$,
we take a specific representative $(U,W,V_\bullet)$ of  the orbit
$\cO_F$ and present a one-parameter family
$(U(\tau),W(\tau),V_\bullet(\tau))$ of subspace configurations
($\tau$ runs over the ground field), such that
$(U(0),W(0),V_\bullet(0))=(U,W,V_\bullet)$, and
$(U(\tau),W(\tau),V_\bullet(\tau))\in\cO_{F'}$ when $\tau\neq 0$.

Since $F'$ is obtained from $F$ by replacing exactly two
indecomposable summands with some other object (consisting of one, two
or three indecomposables), and all the other summands in $F$ remain
unchanged, we can assume that $F$ consists only of these two objects. It
turns out to be convenient to take the representative
$(U,W,V_\bullet)$ in its standard form, as indicated in the
beginning of this subsection.

Now consider all the cases listed in Section~\ref{moves}. We will
consider an initial pair of objects depending on numbers
$i,j,i',j\in[0,n]\cup\{\infty\}$, where $n=\dim V$; when we need to speak about linear
combinations of vectors involving $e_0$ or $e_\infty$, we follow the
convention from the beginning of this subsection. By $V_\bullet$ we
always denote the flag whose components are spanned by basis vectors
$\{e_1,\dots,e_n\}$, such that $\dim V_{a-\alpha}/V_{a_{\alpha-1}}=1$ iff
$\alpha\in\{i,j,i',j'\}$, and 0 otherwise. This flag will always be
invariant along the curves we are going to construct:
$V_\bullet(\tau)=V_\bullet$.

{\bf I.a)} $F=I_{ij}\oplus I_{i'j'}$, $F'=I_{i'j}\oplus I_{ij'}$,
where $i'<i<j'<j$.
$$
(U,W)=(\langle e_j,e_{j'}\rangle,\langle e_i+e_j,e_{i'}+e_{j'}\rangle),
$$
$$
(U(\tau),W(\tau))=(\langle
e_j,e_{j'}\rangle,\langle e_i+e_j,e_{i'}+e_{j'}+\tau e_{j}\rangle).
$$
The triple $(U(\tau),W(\tau),V_\bullet)$ for each nonzero $\tau$
corresponds to the object $F'=I_{i'j}\oplus I_{ij'}$, as may be
seen by calculating its rank numbers, or by
the decomposition of this configuration into a direct sum of two
indecomposables.

Note that this deformation also works for the case when $i'=0$
or/and $j=\infty$.

{\bf I.b)} $F=I_{ij}\oplus I_{i'j'}$, $F'=I_{i'j}\oplus
I_{i}^+\oplus I_{j'}^-$ or $F'=I_{i'j}\oplus I_{i}^-\oplus
I_{j'}^+$ where $i'<j'\leq i<j$. In the first case the initial
configuration
$$
(U,W)=(\langle e_{j'},e_{j}\rangle,\langle
e_{i'}+e_{j'},e_i+e_j\rangle),
$$
is deformed to
$$
(U(\tau),W(\tau))=(\langle e_{j'}+\tau e_i,e_{j}\rangle,\langle
e_{i'}+e_{j'},e_i+e_j\rangle).
$$
and in the second one~--- to
$$
(U(\tau),W(\tau))=(\langle e_{j'},e_{j}\rangle,\langle
e_{i'}+e_{j'}+\tau e_{i},e_i+e_j\rangle).
$$

{\bf I.c)} 
$F=I_{ij}\oplus I^+_{i'}$, $F'=I_{i}^+\oplus I_{i'j}$, where
$i'<i<j$.
$$
(U,W)=(\langle e_{i'},e_{j}\rangle,\langle e_i+e_j\rangle),
$$
$$
(U(\tau),W(\tau))=(\langle e_{i'}+\tau e_i,e_{j}\rangle,\langle
e_i+e_j\rangle).
$$

Similarly, if $F=I_{ij}\oplus I^-_{i'}$ for $i'<i<j$, this object
is transformed to $F'=I_{i}^-\oplus I_{i'j}$: for the
representative
$$
(U,W)=(\langle e_{j}\rangle,\langle e_{i'},e_i+e_j\rangle)
$$
there is a curve
$$
(U(\tau),W(\tau))=(\langle e_{j}\rangle,\langle e_{i'}+\tau e_i,
e_i+e_j\rangle),
$$
having the configuration type $F'$.

{\bf I.d)} $F=I_{i'j'}\oplus I^{+}_{j}$ for $i'<j'<j$, and
$F'=I_{j'}^{+}\oplus I_{i'j}$. Similarly,
$$
(U,W)=(\langle e_{j'},e_j\rangle,\langle e_{i'}+e_{j'}\rangle),
$$
and
$$
(U(\tau),W(\tau))=(\langle e_{j'},e_j\rangle,\langle e_{i'}+e_{j'}+\tau e_j\rangle).
$$

For $F=I_{i'j'}\oplus I^{-}_{j}$ for $i'<j'<j$, and
$F'=I_{j'}^{-}\oplus I_{i'j}$, we have
$$
(U,W)=(\langle e_{j'}\rangle,\langle e_{i'}+e_{j'},e_j\rangle),
$$
$$
(U(\tau),W(\tau))=(\langle e_{j'}+\tau e_j\rangle,\langle e_{i'}+e_{j'}, e_j\rangle).
$$

{\bf I.e)} $F=I^+_{i}\oplus I^-_{i'}$ for $i'<i$,
$F'=I_{i'i}$.
$$
(U,W)=(\langle e_{i}\rangle,\langle e_{i'}\rangle),
$$
and
$$
(U(\tau),W(\tau))=(\langle e_{i}\rangle,\langle e_{i'}+\tau
e_{i}\rangle).
$$
The case $F=I^-_{i}\oplus I^+_{i'}$, $F'=I_{i'i}$ for $i'<i$
is completely analogous.

And here comes the last case:

{\bf II.} $F=I_{ij}\oplus I_{i'j'}$, where $0\leq
i'<j'<i<j\leq\infty$, and $F'=I_{i'i}\oplus I_{j'j}$. Then
$$
(U,W)=(\langle e_{j'},e_{j}\rangle,\langle
e_{i'}+e_{j'},e_i+e_j\rangle),
$$
and
$$
(U(\tau),W(\tau))=(\langle e_{j'}+\tau e_i,e_j\rangle,\langle
e_{i'}+e_{j'}+\tau e_i,e_i+e_j\rangle)
$$

So, for all the possible types of elementary moves we constructed
curves that are contained in the closure of the ``larger"
orbit and that intersect the ``smaller" orbit in exactly one point.
This proves the lemma.
\endproof

\proof[Proof of Prop.~\ref{parab}] Each minimal parabolic subgroup may be presented as the closure of the product
$$
P_i=\overline{U_{i}^-\cdot B},
$$
where $U^-_i=\{E+\tau E_{i+1,i}\mid \tau\in\KK\}$ is a one-dimensional unipotent subgroup consisting of the matrices whose diagonal entries equal 1, and the only nonzero non-diagonal entry, situated in the $i+1$-th line and $i$-th column, equals $\tau$.

For a pair of orbits $\cO_Y$ and $\cO_Y'$, such that $\overline{\cO_Y'}=P_i\overline{\cO_Y}$, and a representative $(U,W)\in \cO_Y$, the action of $U_i^-$ gives us the curve  $U^-_i(U,W)=\{(U(\tau),W(\tau)\}\subset\overline{\cO_Y'}$. For a general $\tau$, the point $(U(\tau),W(\tau))$ belongs to the orbit $\cO'_Y$.

We see that, for the canonical representative $(U,W,V_\bullet)\in\cO\subset X$ corresponding to $\cO_Y\subset Y$, the curve $(U(\tau),W(\tau),V_\bullet)\subset\cO'$ is exactly the one that was constructed in the proof of Lemma~\ref{mvdeg}. The corresponding region has its source and sink on the roads beginning at $I_{i+1,\infty}$ and $I_{i\infty}$ and is \emph{not} of type I.b).

Conversely, let $F\lessdot F'$. Suppose that the elementary move transferring $F$ to $F'$ is not of type I.b), and that the source and the sink of the corresponding minimal admissible region belong to the roads beginning in $I_{r\infty}$ and $I_{s\infty}$ respectively, $s<r$. Then the curve constructed in the proof of Lemma~\ref{mvdeg} is of the form
\begin{eqnarray*}
U(\tau)&=&A_{rs}(\tau)U;\\
W(\tau)&=&A_{rs}(\tau)W;\\
V_\bullet(\tau)&=&V_\bullet,
\end{eqnarray*}
where $A_{rs}(\tau)=E+\tau E_{rs}$ is again a matrix with one nonzero nondiagonal entry.  So, this action is given by the minimal parabolic subgroup $P_i$ iff $s=i$ and $r=i+1$.
%
\endproof

\subsection{Degeneration order implies rank order}\label{ssdegrk}

\proof[Proof of Lemma \ref{degrk}] According to Proposition
\ref{ranks}, it suffices to show that all the inequalities of the
form
\begin{eqnarray}
\dim V_{a_i}\cap U&\geq& d;\nonumber\\
\dim V_{a_i}\cap W&\geq& d;\nonumber\\
\dim V_{a_i}\cap U\cap W&\geq& d;\nonumber\\
\dim (((U\cap V_{a_j})+(W\cap V_{a_j}))\cap V_{a_i})+\dim(U\cap W\cap V_{a_j})&\geq&
d\label{rank_ij}
\end{eqnarray}
\noindent define closed conditions on $X=\Gr(k,V)\times\Gr(l,V)\times\Fl_\ba(V)$.

For the first three families of inequalities this is clear --- these
conditions define closed subvarieties in
$X$ cut out by vanishing of
certain determinants in the homogeneous coordinates on $X$. Let us
show this for the last family of inequalities.

Fix $i$ and $j$, $i<j$, and take a configuration of subspaces
$(U,W,V_\bullet)$. Now define a linear map
$$
\varphi_{ij}\colon
(U\cap V_{a_j})\times(W\cap V_{a_j})\to V_{a_j}/V_{a_i}
$$
by
$$
(u,w)\mapsto u+w\mod V_{a_i}.
$$

The dimension of its kernel equals $\dim (((U\cap V_{a_j})+(W\cap V_{a_j}))\cap V_{a_i})+\dim(U\cap W\cap V_{a_j})$. Indeed,
\begin{multline*}
\dim \Ker(\varphi_{ij})=\dim(U\cap V_{a_j})+
\dim(W\cap V_{a_j})-\rk \varphi_{ij}=\\
\dim(U\cap V_{a_j})+\dim (W\cap V_{a_j})-\dim(((U\cap V_{a_j})+(W\cap V_{a_j}))/V_{a_i})=\\
\dim(U\cap V_{a_j})+\dim (W\cap V_{a_j})-\dim((U\cap V_{a_j})+(W\cap V_{a_j}))+\\\dim(((U\cap V_{a_j})+(W\cap V_{a_j}))\cap V_{a_i})=\\
\dim((U\cap V_{a_j})\cap(W\cap V_{a_j}))+\dim(((U\cap V_{a_j})+(W\cap V_{a_j}))\cap V_{a_i})=\\
\dim(U\cap W\cap V_{a_j})+\dim(((U\cap V_{a_j})+(W\cap V_{a_j}))\cap V_{a_i}).
\end{multline*}

Now let us prove that the condition $\dim\Ker\varphi_{ij}\ge d$ defines a closed condition on $X$. This will be done as follows. Consider the direct product $Y$ of $X$ and three copies of $V=V_n$:
$$
Y=\Gr(k,V)\times\Gr(l,V)\times\Fl_\ba(V)\times V\times V\times V,
$$
and take the subset $Z_{ij}\subset Y$ formed by the sixtuples $(U,W,V_\bullet,x,y,z)\in Y$ satisfying the following conditions:

\begin{eqnarray*}
x,y\in V_{a_j};&\\
x\in U;&\\
y\in W;&\\
z\in V_{a_i};&\\
x+y=z &\text{(as vectors in $V$)}.
\end{eqnarray*}

Obviously, $Z_{ij}$ is closed in $Y$. Moreover, $\Ker\varphi_{ij}\simeq\pi^{-1}_{ij}((U,W,V_\bullet))$, where $\pi_{ij}$ is the projection $Z_{ij}\to X$.

This means that the condition \ref{rank_ij} is equivalent to the condition
$$
\dim\pi^{-1}_{ij}((U,W,V_\bullet))\ge d,
$$
and the latter condition is closed on $X$.
\endproof

\subsection{Rank order implies move order}\label{ssrkmv}

Let us first establish two general facts about rank numbers.

\begin{prop} The set of rank numbers uniquely
defines the corresponding object.
\end{prop}

\proof Assume the contrary: let $F$ and $F'$ correspond to the same
set of rank numbers. This means that $\langle I,F\rangle=\langle
I,F'\rangle$ for each indecomposable $I$.

Since the direct sums of objects correspond to the sums of their
rank numbers, one can consider that no indecomposable objects appear
in $F$ and $F'$ simultaneously. Now take two rightmost objects $I$
and $I'$ (in the sense of AR-quiver of type $D$) occuring in $F$
and $F'$. Without loss of generality suppose that $I$ is situated in
the same column or to the right of $I'$, and, consequently,
(non-strictly) to the right of all indecomposable objects appearing
in $F'$. This means that $\langle I,F'\rangle=0$. Similarly, $I$ is
situated non-strictly to the right of all the indecomposables from
$F$, except for $I$ itself. So $\langle I,F\rangle=\langle
I,I\rangle=1$, a contradiction.\endproof

\begin{prop}\label{cornersums} Let $\fA$ be a region with initial vertices
$I_1$ (source) and $I_2$ (sink), and $J$ the sum of the indecomposable objects corresponding to the terminal vertices of $\fA$. Then for an arbitrary object $F$
$$
\langle I_1,F\rangle+\langle I_2,F\rangle\geq\langle J,F\rangle.
$$
Moreover, if $\fA\setminus I_2$ contains no indecomposable subobject of
$F$, the inequality is an equality.
\end{prop}

\proof By bilinearity of $\langle\cdot,\cdot\rangle$, one can assume $F$ to be indecomposable. So, suppose $F=I$.

Let $I'$ and $I''$ be two neighbor indecomposable objects in a horizontal line (that is, $I_{ij}$ and $I_{i+1,j+1}$, or $I_i^{\pm}$ and $I_{i+1}^{\mp}$). Also denote by $J$ the sum of the objects corresponding to vertices situated on the paths from $I'$ to $I''$ ($J$ may consist of at most three indecomposable objects). With (\ref{expranks}) from Page~\pageref{expranks}, one can see that
\begin{equation}
\langle I',I\rangle+\langle I'',I\rangle\geq \langle J,I\rangle,\label{elem}
\end{equation}
and the inequality is strict iff $I'=I$.

Now, taking the sum of the inequalities (\ref{elem}) over all pairs $(I',I'')$, where both $I'$ and $I''$ belong to $\fA$, we obtain the desired inequality. If all the inequalities (\ref{elem}) are equalities, the latter is equality as well.
\endproof

Next, we need notions of the \emph{interior} and the
\emph{nucleus} of a region. 

\defin Let $\fA$ be a region. The \emph{interior} and the \emph{nucleus} of $\fA$ (denoted by $\Int\fA$ and $\Nuc\fA$, respectively) are sets of indecomposable objects, defined as follows:
\begin{eqnarray*}
\Int\fA&=&\{I\mid \sum_{I'\in\Term\fA}\langle I,I'\rangle<\sum_{I'\in\Init\fA}\langle I,I'\rangle\};\\
\Nuc\fA&=&\{I\mid \sum_{I'\in\Term\fA}\langle I,I'\rangle=\sum_{I'\in\Init\fA}\langle I,I'\rangle-2\}\subset\Int\fA;\\
\end{eqnarray*}
A simple verification shows that $\Int\fA\subset\fA$ and that the difference between $\langle I, F\rangle$ and $\langle I, \tilde F\rangle$ does not exceed 1 for regions of type I.a)-e) and 2 for regions of type II. (So, the nucleus is nonempty only for regions of type II).

On the figures below, for a region of each type its
nucleus is marked with stars, and the interior is formed by the
union of the nucleus with the set of black dots. As before, the initial and terminal vertices are outlined by  squares and circles, respectively.

{\bf I.a)}
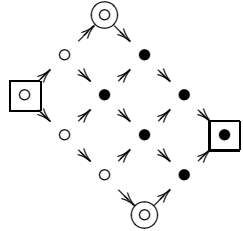
$$
\xymatrix@ur @=6pt{
*+[F]{\circ}\ar[r]\ar[d]&{\circ}\ar[r]\ar[d]&*+[o][F]{\circ}\ar[d]\\
{\circ}\ar[r]\ar[d]&{\bullet}\ar[r]\ar[d]&{\bullet}\ar[d]\\
{\circ}\ar[r]\ar[d]&{\bullet}\ar[r]\ar[d]&{\bullet}\ar[d]\\
*+[o][F]{\circ}\ar[r]&{\bullet}\ar[r]&*+[F]{\bullet} \\
}
$$

{\bf I.b)}
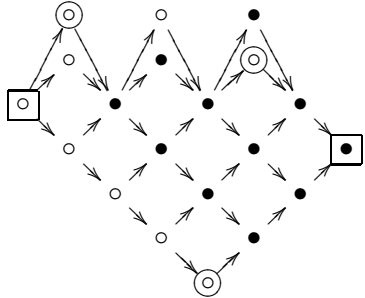
$$
\xymatrix@=6pt{
& *+[o][F]{\circ}\ar[ddr] && {\circ}\ar[ddr] & &{\bullet}\ar[ddr]\\
& {\circ}\ar[dr] &&{\bullet}\ar[dr] && *+[o][F]{\circ}\ar[dr]\\
*+[F]{\circ}\ar[uur] \ar[ur] \ar[dr] && {\bullet} \ar[ur] \ar[dr]\ar[uur] &&
{\bullet} \ar[uur] \ar[ur]\ar[dr]
&& {\bullet}\ar[dr]\\
&{\circ}\ar[ur] \ar[dr] && {\bullet}\ar[ur] \ar[dr] &&
{\bullet}\ar[ur] \ar[dr] && *+[F]{\bullet}\\
&&{\circ}\ar[ur]\ar[dr] && {\bullet}\ar[ur]\ar[dr] && {\bullet} \ar[ur]\\
&&&{\circ}\ar[ur]\ar[dr] && {\bullet}\ar[ur]\\
&&&&*+[o][F]{\circ}\ar[ur] }
$$

{\bf I.c)}
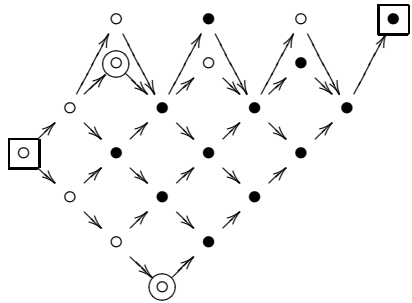
$$
\xymatrix@=6pt{
& &&{\circ}\ar[ddr]&&{\bullet}\ar[ddr] && {\circ}\ar[ddr] &&*+[F]{\bullet}\\
& &&*+[o][F]{\circ}\ar[dr] && {\circ}\ar[dr]&&{\bullet}\ar[dr]\\
 && {\circ} \ar[ur]\ar[uur] \ar[dr] && {\bullet} \ar[uur]\ar[ur] \ar[dr]
&& {\bullet}\ar[uur] \ar[ur] \ar[dr] && {\bullet}\ar[uur]\\
&*+[F]{\circ}\ar[ur] \ar[dr] && {\bullet}\ar[ur] \ar[dr] &&
{\bullet}\ar[ur] \ar[dr]
&& {\bullet}\ar[ur]\\
&&{\circ}\ar[ur]\ar[dr] && {\bullet}\ar[ur]\ar[dr] && {\bullet} \ar[ur]\\
&&&{\circ}\ar[ur]\ar[dr] && {\bullet}\ar[ur]\\
&&&&*+[o][F]{\circ}\ar[ur]}
$$

{\bf I.d)}
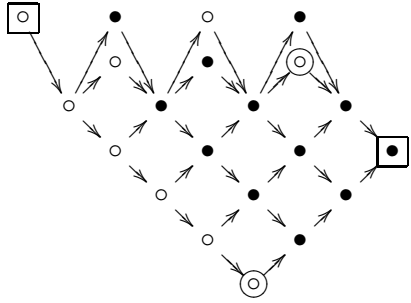
$$
\xymatrix@=6pt{
*+[F]{\circ}\ar[ddr] &&{\bullet}\ar[ddr]&&{\circ}\ar[ddr] && {\bullet}\ar[ddr]\\
&&{\circ}\ar[dr] && {\bullet}\ar[dr]&&*+[o][F]{\circ}\ar[dr]\\
& {\circ} \ar[ur]\ar[uur]\ar[dr] && {\bullet} \ar[uur]\ar[ur] \ar[dr]
&& {\bullet}\ar[uur] \ar[ur] \ar[dr] && {\bullet}\ar[dr]\\
&&{\circ}\ar[ur] \ar[dr] && {\bullet}\ar[ur] \ar[dr] &&
{\bullet}\ar[ur] \ar[dr] && *+[F]{\bullet}\\
&&&{\circ}\ar[ur]\ar[dr] && {\bullet}\ar[ur]\ar[dr] && {\bullet} \ar[ur]\\
&&&&{\circ}\ar[ur]\ar[dr] && {\bullet}\ar[ur]\\
&&&&&*+[o][F]{\circ}\ar[ur]}
$$

{\bf I.e)}
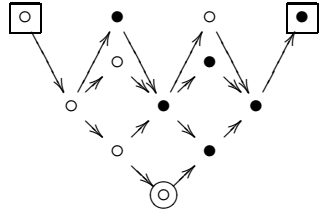
$$
\xymatrix@=6pt{
*+[F]{\circ}\ar[ddr]&&{\bullet}\ar[ddr]&&{\circ}\ar[ddr]&&*+[F]{\bullet}\\
&&{\circ}\ar[dr]&&{\bullet}\ar[dr]\\
&{\circ}\ar[uur]\ar[ur]\ar[dr]&&{\bullet}\ar[uur]\ar[ur]\ar[dr]&&{\bullet}\ar[uur]\\
&&{\circ}\ar[ur]\ar[dr]&&{\bullet}\ar[ur]\\
&&&*+[o][F]{\circ}\ar[ur] }
$$

{\bf II.} 
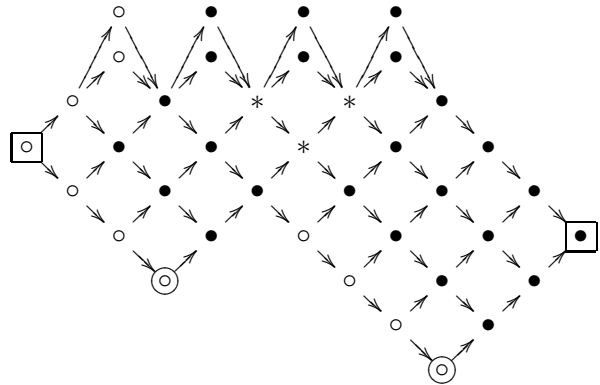
$$
\xymatrix@=6pt{
&& {\circ}\ar[ddr] &&{\bullet}\ar[ddr] &&{\bullet}\ar[ddr]&&{\bullet}\ar[ddr]\\
&& {\circ}\ar[dr] &&{\bullet}\ar[dr] &&{\bullet}\ar[dr]&&{\bullet}\ar[dr]\\
&{\circ}\ar[uur] \ar[ur] \ar[dr] && {\bullet} \ar[uur] \ar[ur]
\ar[dr] && {*} \ar[uur] \ar[ur] \ar[dr]&& {*} \ar[uur] \ar[ur]
\ar[dr]&& {\bullet}\ar[dr]\\
*+[F]{\circ}\ar[ur] \ar[dr] && {\bullet}\ar[ur] \ar[dr] &&
{\bullet}\ar[ur] \ar[dr] &&
{*}\ar[ur] \ar[dr] && {\bullet}\ar[ur] \ar[dr] && {\bullet}\ar[dr]\\
&{\circ}\ar[ur] \ar[dr] && {\bullet}\ar[ur] \ar[dr] &&
{\bullet}\ar[ur] \ar[dr] &&
{\bullet}\ar[ur] \ar[dr] && {\bullet}\ar[ur] \ar[dr] && {\bullet}\ar[dr]\\
&&{\circ}\ar[ur] \ar[dr] && {\bullet}\ar[ur] &&
{\circ}\ar[ur]\ar[dr] &&
{\bullet}\ar[ur] \ar[dr] && {\bullet}\ar[ur] \ar[dr] && *+[F]{\bullet}\\
&&&*+[o][F]{\circ}\ar[ur]&&&&{\circ}\ar[ur]\ar[dr]&&{\bullet}\ar[ur]\ar[dr]&&{\bullet}\ar[ur]\\
&&&&&&&&{\circ}\ar[ur]\ar[dr]&&{\bullet}\ar[ur]\\
&&&&&&&&&*+[o][F]{\circ}\ar[ur] }
$$

Now let us pass to the proof of Lemma \ref{rkmv}.

\proof[Proof of Lemma \ref{rkmv}] Let  $F$ and $F'$ be two objects,
such that $\underline\dim F=\underline\dim F'$ and
$F\rkless F'$. We have
$\langle I,F\rangle\geq\langle I,F'\rangle $ for all indecomposables $I$. For the ``fake vertex" $I_{0\infty}$ we set
$\langle I_{0\infty},F\rangle=\langle I_{0\infty},F'\rangle =0$.

We begin with the following definition, which will be the last one in this paper.

\defin A region $\fB$ is said to be \emph{dominant} w.r.t. $F$ and $F'$, if
the following inequalities hold:
$$
\langle I,F\rangle>\langle I,F'\rangle \qquad \forall I\in\Int\fB;
$$
$$
\langle I,F\rangle>\langle I,F'\rangle+1 \qquad \forall I\in\Nuc\fB.
$$
(Of course, the second set of inequalities is trivial for regions of type I).


The following technical lemma is essential for the sequel.

\begin{lemma}\label{tech} With the notation as above, 
take a rightmost object $I$, such that the corresponding rank
numbers for $F$ and $F'$ differ: $\langle I,F\rangle>\langle
I,F'\rangle$. Then there exists a dominant region $\fB$ with
 sink $I$ and an indecomposable object $J\neq I$ situated in $\fB$
and occuring in $F$ as a direct summand.
\end{lemma}

\proof Take a maximal dominant region $\fB$ with  sink $I$. Assume the
contrary: no indecomposable summand of $F$ other than $I$ is situated
in $\fB$.

1. First suppose that $\fB$ is of type II, with sink $I=I_{ij}$ and source $I'=I_{i'j'}$. We know that $i<j<i'<j'$.

Since $\fB$ is maximal dominant, there must exist two objects $J_1$ and $J_2$ with the property
$$
\langle J_{1,2},F\rangle=\langle J_{1,2},F'\rangle,
$$
such that 
$$
J_1\in\{I_{\alpha j'}\mid \alpha\in [j,i')\}
$$
and 
$$
J_2\in\{I_{\beta i'}\mid \beta\in (i,j]\}\cup \{I_{i'\gamma}\mid\gamma\in (i',j')\}\cup\{I^{\pm}_{i'}\}
$$
(otherwise $\fB$ would be contained in a larger dominant region).

According to the position of $J_2$, three cases can occur:

\medskip

1a. $J_1=I_{\alpha j'}$, $J_2=I_{\beta i'}$, where $\alpha\in [j,i')$, $\beta\in (i,j]$. 

Consider also two objects $I_{i'j'}$ and $I_{\beta\alpha}$. These four objects determine a region of type II:
$$
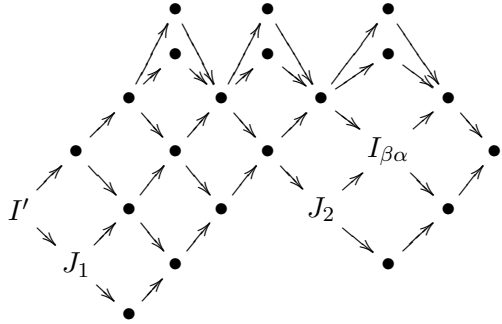
\xymatrix@=6pt{
&&&\bullet\ar[ddr] &&\bullet\ar[ddr] &&\bullet\ar[ddr]\\
&&&\bullet\ar[dr] &&\bullet\ar[dr] &&\bullet\ar[dr]\\
&&\bullet\ar[ur]\ar[uur]\ar[dr] &&\bullet\ar[ur]\ar[uur]\ar[dr] &&\bullet\ar[dr]\ar[ur]\ar[uur] &&\bullet\ar[dr]\\
&\bullet\ar[ur]\ar[dr] &&\bullet\ar[ur]\ar[dr] &&\bullet\ar[dr]\ar[ur] &&{I_{\beta\alpha}}\ar[ur]\ar[dr] &&\bullet\\
I'\ar[ur]\ar[dr] &&\bullet\ar[ur]\ar[dr] &&\bullet\ar[ur] && J_2\ar[ur]\ar[dr] &&\bullet\ar[ur]\\
& J_1\ar[ur]\ar[dr] && \bullet\ar[ur] &&&&\bullet\ar[ur]\\
&&\bullet\ar[ur]
}
$$
Apply Prop.~\ref{cornersums} twice to this region, taking into account that $I_{\beta\alpha}\in\Int\fB$:
\begin{multline*}
\langle I_{i'j'},F\rangle=\langle J_1,F\rangle+\langle J_2,F\rangle-\langle I_{\beta\alpha},F\rangle\\
<\langle J_1,F'\rangle+\langle J_2,F'\rangle-\langle I_{\beta\alpha},F'\rangle\leq \langle I_{i'j'}, F'\rangle,
\end{multline*}
that gives us a contradiction. This means that this smaller region, and hence $\fB$, contain subobjects of $F$ different from $I$.

\medskip

1b. $J_1=I_{\alpha j'}$, $J_2=I_{i'\gamma}$, where $\alpha\in [j,i')$, $\gamma\in (i',j')$.

In this case, we consider the objects $I_{i'j'}$ and $I_{\alpha\gamma}$:
$$
\xymatrix@=6pt{
&&&\bullet\ar[ddr] &&\bullet\ar[ddr] &&\bullet\ar[ddr]\\
&&&\bullet\ar[dr] &&\bullet\ar[dr] &&\bullet\ar[dr]\\
&&J_2\ar[ur]\ar[uur]\ar[dr] &&\bullet\ar[ur]\ar[uur]\ar[dr] &&\bullet\ar[dr]\ar[ur]\ar[uur] &&\bullet\ar[dr]\\
&\bullet\ar[ur]\ar[dr] &&I_{\alpha\gamma}\ar[ur]\ar[dr] &&\bullet\ar[dr]\ar[ur] &&\bullet\ar[ur]\ar[dr] &&\bullet\\
I'\ar[ur]\ar[dr] &&\bullet\ar[ur]\ar[dr] &&\bullet\ar[ur] && \bullet\ar[ur]\ar[dr] &&\bullet\ar[ur]\\
& J_1\ar[ur]\ar[dr] && \bullet\ar[ur] &&&&\bullet\ar[ur]\\
&&\bullet\ar[ur]
}
$$
 and again apply the same Proposition:
\begin{multline*}
\langle I_{i'j'},F\rangle=\langle J_1,F\rangle+\langle J_2,F\rangle-\langle I_{\alpha\gamma},F\rangle\\
<\langle J_1,F'\rangle+\langle J_2,F'\rangle-\langle I_{\alpha\gamma},F'\rangle\leq \langle I_{i'j'}, F'\rangle,
\end{multline*}
obtaining a contradiction with our assumption.

\medskip

1c. $J_1=I_{\alpha j'}$, $\alpha\in [j,i')$, and $J_2=I^{\pm}_{i'}$. 

We consider the pair of objects $(I_{i'j'},I^{\pm}_\alpha)$ and again apply the same procedure (see figure below).
$$
\xymatrix@=6pt{
&&&J_2\ar[ddr] &&\bullet\ar[ddr] &&\bullet\ar[ddr]\\
&&&\bullet\ar[dr] &&I_\alpha^{\pm}\ar[dr] &&\bullet\ar[dr]\\
&&\bullet\ar[ur]\ar[uur]\ar[dr] &&\bullet\ar[ur]\ar[uur]\ar[dr] &&\bullet\ar[dr]\ar[ur]\ar[uur] &&\bullet\ar[dr]\\
&\bullet\ar[ur]\ar[dr] &&\bullet\ar[ur]\ar[dr] &&\bullet\ar[dr]\ar[ur] &&\bullet\ar[ur]\ar[dr] &&\bullet\\
I'\ar[ur]\ar[dr] &&\bullet\ar[ur]\ar[dr] &&\bullet\ar[ur] && \bullet\ar[ur]\ar[dr] &&\bullet\ar[ur]\\
& J_1\ar[ur]\ar[dr] && \bullet\ar[ur] &&&&\bullet\ar[ur]\\
&&\bullet\ar[ur]
}
$$

\bigskip

2. The region $\fB$ is of type I.a)--I.c). This means that its source $I'$ is of the form $I_{ij}$. 

The maximality of $\fB$ implies the existence of at least two objects $J\in\fB$, such that $\langle J,F\rangle=\langle J,F'\rangle$. We distinguish between the following subcases:



\medskip

2a. There are two such objects of the form $J_1=I_{i'j}$ and $J_2=I_{ij'}$, $j'\in(i,j)$. Then we can consider the objects $I_{ij}$ and $I_{i'j'}$:
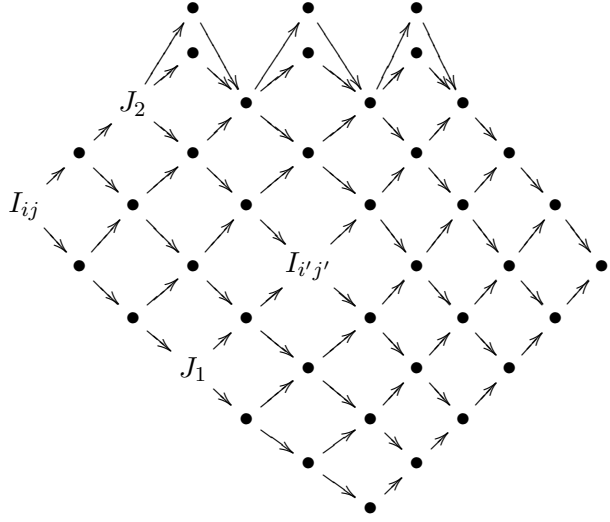
$$
\xymatrix@=6pt{
&&&\bullet\ar[ddr] &&\bullet\ar[ddr]&&\bullet\ar[ddr] \\
&&&\bullet\ar[dr] && \bullet\ar[dr]&&\bullet\ar[dr]\\
&&J_2 \ar[ur]\ar[uur]\ar[dr] &&\bullet\ar[uur]\ar[ur]\ar[dr] &&\bullet\ar[dr]\ar[ur]\ar[uur] &&\bullet\ar[dr]\\
&\bullet\ar[ur]\ar[dr] &&\bullet\ar[ur]\ar[dr] &&\bullet\ar[dr]\ar[ur] &&\bullet\ar[ur]\ar[dr] &&\bullet\ar[dr]\\
I_{ij}\ar[ur]\ar[dr] &&\bullet\ar[ur]\ar[dr] &&\bullet\ar[dr]\ar[ur] &&\bullet\ar[ur]\ar[dr] &&\bullet\ar[dr]\ar[ur] &&\bullet\ar[dr]\\
&\bullet\ar[ur]\ar[dr] &&\bullet\ar[ur]\ar[dr] &&I_{i'j'}\ar[dr]\ar[ur] &&\bullet\ar[ur]\ar[dr] &&\bullet\ar[dr]\ar[ur] &&\bullet\\
&&\bullet\ar[ur]\ar[dr] &&\bullet\ar[ur]\ar[dr] &&\bullet\ar[dr]\ar[ur] &&\bullet\ar[ur]\ar[dr] &&\bullet\ar[ur]\\
&&&J_1\ar[ur]\ar[dr] &&\bullet\ar[ur]\ar[dr] &&\bullet\ar[dr]\ar[ur]  &&\bullet\ar[ur]\\
&&&&\bullet\ar[ur]\ar[dr] &&\bullet\ar[ur]\ar[dr] &&\bullet\ar[ur]\\
&&&&&\bullet\ar[ur]\ar[dr] &&\bullet\ar[ur]\\
&&&&&&\bullet\ar[ur]\\
}
$$
and apply Prop.~\ref{cornersums} twice, writing
\begin{multline*}
\langle I_{ij},F\rangle=\langle J_1,F\rangle+\langle J_2,F\rangle-\langle I_{i'j'},F\rangle\\
<\langle J_1,F'\rangle+\langle J_2,F'\rangle-\langle I_{i'j'},F'\rangle\leq \langle I_{ij}, F'\rangle.
\end{multline*}
This gives us a contradiction.
\medskip

2b. $J_1=I_{i'j}$, but for all vertices $I_{ij'}$, where $i<j'<j$, the inequality 
$$
\langle I_{ij'},F\rangle>\langle I_{ij'},F'\rangle
$$
holds. Then, by maximality of $\fB$, there exist two vertices $J_2=I^\pm_i$ and $J_3=I^\mp_{i''}$ (with different signs), such that $I^\pm_i\in\Term\fB$, and
\begin{eqnarray*}
\langle J_2,F\rangle=\langle J_2,F'\rangle\\
\langle J_3,F\rangle=\langle J_3,F'\rangle\\
\end{eqnarray*}
Let us take for $J_3$ the leftmost element of form $I^\mp_\bullet$ situated in $\fB$ and satisfying the latter equality. 

If $i''\leq i'$, we can consider region $\fC$ of type I.c) with $\Init\fC = \{I_{ij}, I^\pm_{i'}\}$ and $\Term\fC=\{I_{i'j}, I^\pm_i\}$, see figure:
$$
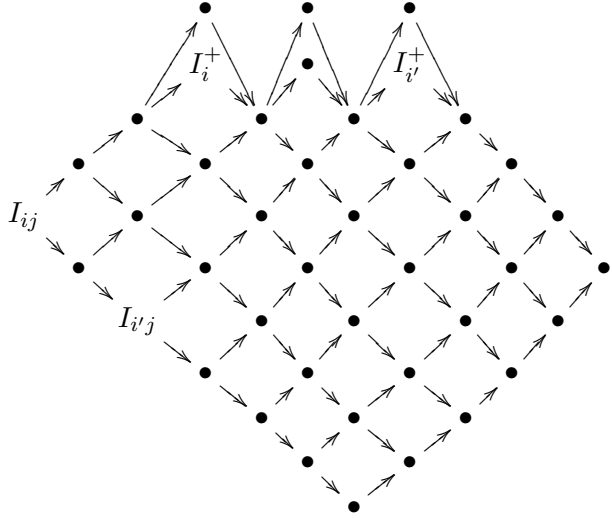
\xymatrix@=6pt{
&&&\bullet\ar[ddr] &&\bullet\ar[ddr]&&{\bullet}\ar[ddr]\\
&&&{I^+_i}\ar[dr] && {\bullet}\ar[dr] &&{I^+_{i'}}\ar[dr]\\
&&\bullet \ar[ur]\ar[uur]\ar[dr] &&\bullet\ar[uur]\ar[ur]\ar[dr] &&\bullet\ar[dr]\ar[ur]\ar[uur] &&\bullet\ar[dr]\\
&\bullet\ar[ur]\ar[dr] &&\bullet\ar[ur]\ar[dr] &&\bullet\ar[dr]\ar[ur] &&\bullet\ar[ur]\ar[dr] &&\bullet\ar[dr]\\
I_{ij}\ar[ur]\ar[dr] &&\bullet\ar[ur]\ar[dr] &&\bullet\ar[dr]\ar[ur] &&\bullet\ar[ur]\ar[dr] &&\bullet\ar[dr]\ar[ur] &&\bullet\ar[dr]\\
&\bullet\ar[ur]\ar[dr] &&\bullet\ar[ur]\ar[dr] &&\bullet\ar[dr]\ar[ur] &&\bullet\ar[ur]\ar[dr] &&\bullet\ar[dr]\ar[ur] &&\bullet\\
&&{I_{i'j}}\ar[ur]\ar[dr] &&\bullet\ar[ur]\ar[dr] &&\bullet\ar[dr]\ar[ur] &&\bullet\ar[ur]\ar[dr] &&\bullet\ar[ur]\\
&&&\bullet\ar[ur]\ar[dr] &&\bullet\ar[ur]\ar[dr] &&\bullet\ar[dr]\ar[ur]  &&\bullet\ar[ur]\\
&&&&\bullet\ar[ur]\ar[dr] &&\bullet\ar[ur]\ar[dr] &&\bullet\ar[ur]\\
&&&&&\bullet\ar[ur]\ar[dr] &&\bullet\ar[ur]\\
&&&&&&\bullet\ar[ur]\\
}
$$
Then we can again apply Prop.~\ref{cornersums} and obtain
\begin{multline*}
\langle I_{ij},F\rangle=\langle I_{i'j},F\rangle+\langle I^\pm_i,F\rangle-\langle I^+_{i'},F\rangle \\ <\langle I_{i'j},F'\rangle+\langle I^\pm_i,F'\rangle-\langle I^+_{i'},F'\rangle\leq \langle I_{ij},F'\rangle.
\end{multline*}

2c. If $i''> i'$, we consider the region $\fC'$ of type I.b), with $\Init\fC'=\{I_{ij}, I_{i'i''}\}$ and $\Term\fC'=\{I_{i'j},I^\pm_i,I^\mp_{i''}\}$,
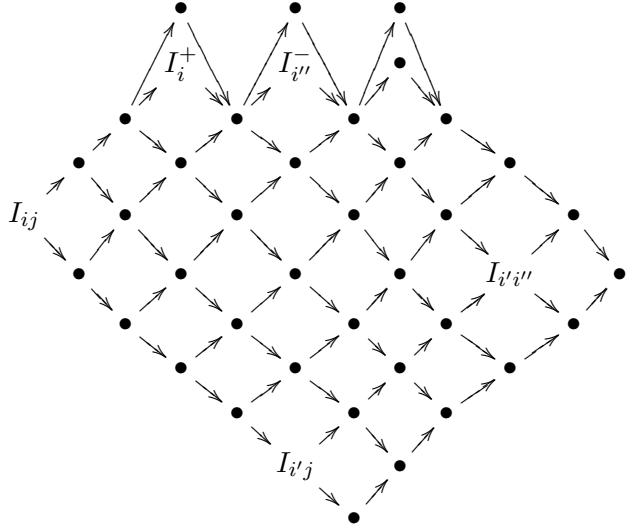
$$
\xymatrix@=6pt{
&&&\bullet\ar[ddr] &&\bullet\ar[ddr]&&{\bullet}\ar[ddr]\\
&&&{I^+_i}\ar[dr] && {I^-_{i''}}\ar[dr] &&{\bullet}\ar[dr]\\
&&\bullet \ar[ur]\ar[uur]\ar[dr] &&\bullet\ar[uur]\ar[ur]\ar[dr] &&\bullet\ar[dr]\ar[ur]\ar[uur] &&\bullet\ar[dr]\\
&\bullet\ar[ur]\ar[dr] &&\bullet\ar[ur]\ar[dr] &&\bullet\ar[dr]\ar[ur] &&\bullet\ar[ur]\ar[dr] &&\bullet\ar[dr]\\
I_{ij}\ar[ur]\ar[dr] &&\bullet\ar[ur]\ar[dr] &&\bullet\ar[dr]\ar[ur] &&\bullet\ar[ur]\ar[dr] &&\bullet\ar[dr]\ar[ur] &&\bullet\ar[dr]\\
&\bullet\ar[ur]\ar[dr] &&\bullet\ar[ur]\ar[dr] &&\bullet\ar[dr]\ar[ur] &&\bullet\ar[ur]\ar[dr] &&{I_{i'i''}}\ar[dr]\ar[ur] &&\bullet\\
&&{\bullet}\ar[ur]\ar[dr] &&\bullet\ar[ur]\ar[dr] &&\bullet\ar[dr]\ar[ur] &&\bullet\ar[ur]\ar[dr] &&\bullet\ar[ur]\\
&&&\bullet\ar[ur]\ar[dr] &&\bullet\ar[ur]\ar[dr] &&\bullet\ar[dr]\ar[ur]  &&\bullet\ar[ur]\\
&&&&\bullet\ar[ur]\ar[dr] &&\bullet\ar[ur]\ar[dr] &&\bullet\ar[ur]\\
&&&&&{I_{i'j}}\ar[ur]\ar[dr] &&\bullet\ar[ur]\\
&&&&&&\bullet\ar[ur]\\
}
$$
Again we apply Prop.~\ref{cornersums} to this region twice, obtaining
\begin{multline*}
\langle I_{ij},F\rangle=\langle I_{i'j},F\rangle+\langle I^\pm_i,F\rangle+\langle I^\mp_{i''},F\rangle-\langle I_{i'i''},F\rangle \\
<\langle I_{i'j},F'\rangle+\langle I^\pm_i,F'\rangle+\langle I^\mp_{i''},F'\rangle-\langle I_{i'i''},F'\rangle <\langle I_{ij},F'\rangle.
\end{multline*}

3. The region $\fB$ is of type I.d) or I.e). Let its source be situated at the vertex $I=I^\pm_j$. The maximality of $\fB$ means that there exists at least one element $I_{ij}$, such that $\langle I_{ij},F\rangle=\langle I_{ij},F'\rangle$. Let $I_{ij}$ be the leftmost element with this property. We distinguish between the following subcases:

3a. There exists an element $I^\pm_{j'}$, such that $\langle I^\pm_{j'},F\rangle=\langle I^\pm_{j'},F'\rangle$, and $i\leq j'$. In this case, take a leftmost such element and consider region $\fC$ of type I.d), defined by $\Init\fC=\{I^\pm_j,I_{ij'}\}$ and $\Term\fC=\{I^\pm_{j'},I_{ij}\}$. It does not contain objects occuring in $F$, so proceed as usual:
\begin{multline*}
\langle I^\pm_{j},F\rangle=\langle I_{ij'},F\rangle-\langle I^\pm_{j'},F\rangle-\langle I_{ij},F\rangle\\
<\langle I_{ij'},F'\rangle-\langle I^\pm_{j'},F'\rangle-\langle I_{ij},F'\rangle\leq \langle I^\pm_{j},F'\rangle,
\end{multline*}
a contradiction.
$$
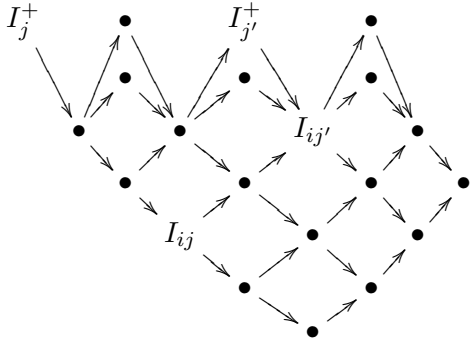
\xymatrix@=6pt{
{I^+_j}\ar[ddr] &&{\bullet}\ar[ddr]&&{I_{j'}^+}\ar[ddr] && {\bullet}\ar[ddr]\\
&&{\bullet}\ar[dr] && {\bullet}\ar[dr]&&{\bullet}\ar[dr]\\
& {\bullet} \ar[ur]\ar[uur]\ar[dr] && {\bullet} \ar[uur]\ar[ur] \ar[dr]
&& {I_{ij'}}\ar[uur] \ar[ur] \ar[dr] && {\bullet}\ar[dr]\\
&&{\bullet}\ar[ur] \ar[dr] && {\bullet}\ar[ur] \ar[dr] &&
{\bullet}\ar[ur] \ar[dr] && {\bullet}\\
&&&{I_{ij}}\ar[ur]\ar[dr] && {\bullet}\ar[ur]\ar[dr] && {\bullet} \ar[ur]\\
&&&&{\bullet}\ar[ur]\ar[dr] && {\bullet}\ar[ur]\\
&&&&&{\bullet}\ar[ur]}
$$

3b. For all elements $I^\pm_{j'}$, such that $i<j'<j$, the inequality $\langle I^\pm_{j'},F\rangle\geq \langle I^\pm_{j'},F'\rangle$ is strict, and the element $I^\mp_i$ belongs to $\fB$. Then we consider $\fC$ of type I.e), with $\Init\fC=\{I_i^\mp,I^\pm_j\}$ and $\Term\fC=\{I_{ij}\}$, and apply the same method:
$$
\langle I^\pm_{j},F\rangle=\langle I^\mp_{i},F\rangle-\langle I_{ij},F\rangle
<\langle I^\mp_{i},F'\rangle-\langle I_{ij},F'\rangle\leq \langle I^\pm_{j},F'\rangle.
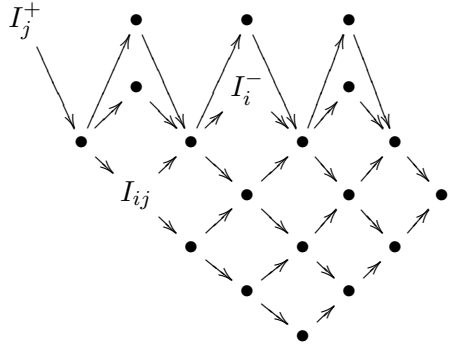
$$
$$
\xymatrix@=6pt{
{I^+_j}\ar[ddr] &&{\bullet}\ar[ddr]&&{\bullet}\ar[ddr] && {\bullet}\ar[ddr]\\
&&{\bullet}\ar[dr] && {I_{i}^-}\ar[dr]&&{\bullet}\ar[dr]\\
& {\bullet} \ar[ur]\ar[uur]\ar[dr] && {\bullet} \ar[uur]\ar[ur] \ar[dr]
&& {\bullet}\ar[uur] \ar[ur] \ar[dr] && {\bullet}\ar[dr]\\
&&{I_{ij}}\ar[ur] \ar[dr] && {\bullet}\ar[ur] \ar[dr] &&
{\bullet}\ar[ur] \ar[dr] && {\bullet}\\
&&&{\bullet}\ar[ur]\ar[dr] && {\bullet}\ar[ur]\ar[dr] && {\bullet} \ar[ur]\\
&&&&{\bullet}\ar[ur]\ar[dr] && {\bullet}\ar[ur]\\
&&&&&{\bullet}\ar[ur]}
$$

3c. Here comes the last possibility: the equality of rank numbers holds in $I_{ij}$, but for all vertices $I^\pm_\alpha\in\fB$, $\alpha\neq j$, the inequality
$$
\langle I^\pm_\alpha, F\rangle\geq \langle I^\pm_\alpha, F'\rangle
$$ 
is strict, and the vertex $I_{i}^\pm$ does not belong to $\fB$. The latter means that $\fB$ is of type I.d) (not I.e)). Denote its sink by $I_{i_0j_0}$.

In this case, we claim that region $\fC$ with $\Init\fC=\{I_{j,j+1},I_{i_0j_0}\}$ and $\Term\fC=\{I_{i_0j}, I_{j+1}\}$ is dominant.

Since $\fB$ is dominant and by the hypothesis of Case 3c, we see that for each $\tilde I\in\Int\fC$, $\langle I,F\rangle\geq\langle I,F'\rangle+1$.

So, we have to show that for each vertex from $\Nuc\fC$, that is, for each vertex of the form $I_{\alpha\beta}$, where $j_0\le\alpha<\beta\leq j-1$, the inequality 
$$
\langle I_{\alpha\beta}, F\rangle\geq \langle I_{\alpha\beta}, F'\rangle+1
$$
is strict.

Let us prove this. Suppose that there exists an object $I_{\alpha_0\beta_0}$, where this
inequality is an equality. Then we can apply
Prop.~\ref{cornersums}, in a slightly different way than before:
\begin{multline*}
\langle I_{\alpha_0j},F\rangle=\langle
I_{\alpha_0\beta_0},F\rangle+\langle I_{\alpha_0j},F\rangle -\langle
I_{i'\beta_0},F\rangle\\<\langle
I_{\alpha_0\beta_0},F'\rangle+1+\langle I_{\alpha_0j},F'\rangle -\langle
I_{i'\beta_0},F'\rangle\leq \langle I_{\alpha_0 j},F'\rangle+1,
\end{multline*}
that yields a contradiction.

Here is the corresponding figure:
$$
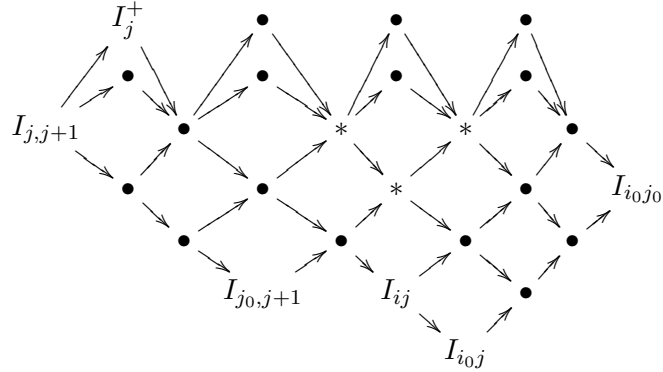
\xymatrix@=6pt{
&{I^+_j}\ar[ddr] &&{\bullet}\ar[ddr]&&{\bullet}\ar[ddr] && {\bullet}\ar[ddr]\\
&{\bullet}\ar[dr]&&{\bullet}\ar[dr] && {\bullet}\ar[dr]&&{\bullet}\ar[dr]\\
{I_{j,j+1}}\ar[ur]\ar[uur]\ar[dr]&& {\bullet} \ar[ur]\ar[uur]\ar[dr] && {*} \ar[uur]\ar[ur] \ar[dr] && {*}\ar[uur] \ar[ur] \ar[dr] && {\bullet}\ar[dr]\\
&{\bullet}\ar[ur]\ar[dr]&&{\bullet}\ar[ur] \ar[dr] && {*}\ar[ur] \ar[dr] &&
{\bullet}\ar[ur] \ar[dr] && {I_{i_0j_0}}\\
&&{\bullet}\ar[ur]\ar[dr]&&{\bullet}\ar[ur]\ar[dr] && {\bullet}\ar[ur]\ar[dr] && {\bullet} \ar[ur]\\
&&&{I_{j_0,j+1}}\ar[ur]&&{I_{ij}}\ar[ur]\ar[dr] && {\bullet}\ar[ur]\\
&&&&&&{I_{i_0j}}\ar[ur]}
$$

So, having obtained a dominant region of type II, we proceed as in the case 1.

The lemma is proved.\endproof

Having such a region $\fB$, let us take a \emph{minimal} dominant region in it; that is, a dominant region $\fC$ satisfying the following properties:

\begin{enumerate}
\item The sink of $\fC$ equals $I$, and its source occurs in $F$ as a direct summand;


\item $\fC$ contains no subobjects of $F$ other that its source and its sink (minimality).
\end{enumerate}

The properties 1 and 2 imply that such a region $\fC$ is minimal admissible. So we may perform the elementary move corresponding to 
$\fC$, thus obtaining an object $\tilde F$ from $F$.
The property of $\fC$ to be dominant implies that $\langle I,\tilde F\rangle\geq\langle I,F'\rangle$ for
each indecomposable object $I$. So, we have found the desired object
$\tilde F$, such that
$$
F\lessdot\tilde F\rkless F'.
$$
This concludes the proof of Lemma \ref{rkmv}.\endproof

\bigskip

{\noindent\footnotesize{\sc Independent University of Moscow, Bolshoi Vlasievskii per., 11, 119002 Moscow, Russia\\
Institut Fourier, 100 rue des Maths, 38400 Saint-Martin d'H\`eres,
France}\\
\emph{E-mail address:} \verb"smirnoff@mccme.ru" }

\end{document}